\documentclass[11pt]{article}   
\usepackage{amssymb,amscd,latexsym}   
\usepackage{amsmath}
\usepackage{amsthm}
\newcommand{\rar}{\rightarrow}
\newcommand{\lar}{\longrightarrow}
\newcommand{\llar}{-\kern-5pt-\kern-5pt\longrightarrow}

\newcommand{\injects}{\hookrightarrow}

\newtheorem{Theorem}{Theorem}[section]
\newtheorem{Lemma}[Theorem]{Lemma}
\newtheorem{Corollary}[Theorem]{Corollary}
\newtheorem{Proposition}[Theorem]{Proposition}
\newtheorem{Remark}[Theorem]{Remark}
\newtheorem{Example}[Theorem]{Example}

\newtheorem{Definition}[Theorem]{Definition}
\newtheorem{Question}[Theorem]{Question}
\def\sqr#1#2{{\vcenter{\hrule height.#2pt
        \hbox{\vrule width.#2pt height#1pt \kern#1pt
            \vrule width.#2pt}
        \hrule height.#2pt}}}
\def\phi{\varphi}
\def\demo{\noindent{\bf Proof. }}
\def\square{\mathchoice\sqr64\sqr64\sqr{4}3\sqr{3}3}
\def\qed{\hspace*{\fill} $\square$}

\DeclareMathOperator{\Hom}{Hom}

\DeclareMathOperator{\beg}{indeg}
\DeclareMathOperator{\Deg}{deg}
\DeclareMathOperator{\MAX}{max}
\DeclareMathOperator{\MIN}{min}

\DeclareMathOperator{\Ext}{Ext}

\DeclareMathOperator{\rank}{rank}

\DeclareMathOperator{\Reg}{reg}

\DeclareMathOperator{\grade}{grade}

\DeclareMathOperator{\END}{end}

\DeclareMathOperator{\E}{\rm E}

\def\XX{{\bf X}}

\def\fm{{\mathfrak m}}

\def\hht{{\rm ht}\,}

\def\edim{{\rm edim}\,}

\def\grade{{\rm grade}\,}
\def\rk{{\rm rank}\,}

\def\restr{{\kern-1pt\restriction\kern-1pt}}

\def\NN{\mathbb N}

\def\pp{{\mathbb P}}

 2

 3

 4

 5

 4

 3

 2

 1

 1

 2

 3

 1

 0

 1

 2

 3

 4

 5

\begin{document}
\begin{center}
{\Large{\bf\sc Plane Cremona maps:}}\\[3pt]
{\Large{\bf\sc saturation and regularity of the base ideal}}\footnotetext{Mathematics Subject Classification 2000
 (MSC2000). Primary  13H15, 13D02, 14E05, 14E07, 14M05;
Secondary 13A02, 13A18, 14B15, 13B22, 13C14, 13D45, 13F30, 14C20.}

\vspace{0.3in}

\hspace{-15pt}{\large\sc Seyed Hamid Hassanzadeh}\footnote{On a  Post-Doc Fellowship (CNPq, Brazil).}
\footnote{This author wishes to thank the hospitality and facilities provided by The Department of Mathematics at
the Universidade Federal de Pernambuco and by IMPA during the preparation of this work.}
\quad\quad\quad
 {\large\sc Aron  Simis}\footnote{Partially
supported by a CNPq grant.}

\end{center}


\bigskip

\begin{abstract}

One studies plane Cremona maps by focusing on the ideal theoretic and homological properties of its
homogeneous base ideal (``indeterminacy locus'').
The {\em leitmotiv} driving a good deal of the work is the relation between the base ideal and its saturation.
As a preliminary one deals with the homological features of arbitrary codimension $2$ homogeneous ideals in a polynomial ring
in three variables over a field which are generated by three forms of the same degree.
The results become sharp when the saturation
is not generated in low degrees, a condition to be given a precise meaning.
An implicit goal, illustrated in low degrees, is a homological classification of  plane Cremona maps according to the
respective homaloidal types.
An additional piece of this work  relates the base ideal of a rational map
to a few additional homogeneous ``companion'' ideals, such as the integral closure, the
$\boldsymbol\mu$-fat ideal and a seemingly novel ideal defined in terms of valuations.

\end{abstract}

\section*{Introduction}

Let $k$ be an  algebraically closed field and let
$R=k[\XX]=k[X_0,\ldots,X_n]$ denote a polynomial ring over $k$,
endowed with the ordinary standard grading.
A rational map $F\colon\pp^n\dasharrow \pp^n$ is defined by
$n+1$ forms in $R$ of the same degree.
If $F$ is birational then it is  called a Cremona map.
Cremona maps are a classical subject
that can do away with any general introduction.
Yet, perhaps less known is a fairly recent  body of results
on the nature and structure of an individual such map -- rather
than on the structure of the Cremona group -- that draws on modern
geometric and algebraic  tools (see, e.g., \cite{CRS},\cite{Pan}, \cite{PanRusso},
\cite{bir2003}, \cite{SiVi}, \cite{SimisVilla}, \cite{{CremonaMexico}}).

In many aspects the properties of the base locus of $F$ play a fundamental role,
where the base locus is the scheme defined by the {\em base ideal}
$I\subset R$ generated by the $n+1$ forms defining $F$.
There is a neat difference between the base ideal and its
ideal theoretic saturation $I^{\rm sat}$. Of course both define the same scheme,
but while $I^{\rm sat}$ serves well the geometric purpose, it is $I$ that gives the nature
of the linear system defining the map.
For plane Cremona maps of degree at least $2$,
 asking when the base ideal $I\subset R$ is saturated is tantamount to asking when
 $R/I$ is a Cohen--Macaulay ring, i.e., when $I$ is generated by the
maximal minors of a $3\times 2$ homogeneous matrix with entries in $R$.
While there are many rational maps on $\pp^2$, in any degree, defined by such ideals, the question for Cremona
maps becomes much tighter.

In this work we deal only with the plane case ($n=2$),
largely focusing on a comparison of these two ideal theoretic versions
of the base scheme.
Furthermore, we exploit a nice interplay between the algebraic properties of the base ideal and the
geometry related to the underlying linear system spanned by the coordinate forms defining the rational map.
Classically, the role of the underlying linear system comes through the  {\em weighted cluster}
associated to it. Exploiting the ideal theoretic side of this linear system we introduce into the
picture other ideals squeezed between the base ideal and the so-called fat ideal associated to the proper points of the
underlying cluster.

\smallskip

Our main results are Theorem~\ref{from_lower_bound}, Theorem~\ref{Rees_of_Jonq}, Theorem~\ref{tcharactrization}
and Theorem~\ref{proper_types_degree5}.

We now proceed to a more detailed description of the sections.

Section~\ref{first} gives a few  homological particulars of a codimension $2$ homogeneous ideal $I\subset R$ generated by
three forms of the same degree $d\geq 2$ in a polynomial ring $R=k[x,y,z]$ over a field $k$.
The use of local cohomology and spectral sequences at this early stage is
justified by a quick derivation of bounds on the regularity of $R/I$ in terms of $d$ and some sharp prediction for the
Betti numbers of the corresponding minimal free resolution.
Furthermore, the true impact of these preliminaries is toward the relation between $I$ and $I^{\rm sat}$.
Here a crucial assumption is that the initial degree of the $R$-module $I^{\rm sat}/I$ is at least $d+1$.
Coupled with the homological preliminaries this hypothesis triggers further lower bounds for $d$ if $I$ is non-saturated,
gives a sharp upper bound for the saturation exponent in terms of the regularity,
and allows to establish the Betti numbers of the resolution of $R/I$ for $d\leq 7$ (see Theorem~\ref{from_lower_bound} and
Proposition~\ref{P567}).

 We introduce a couple of ideals naturally related to the base ideal $I$ of a plane rational map.
One of these is a version of the usual fat ideal for the case when the rational map has infinitely near base points.
Of course, infinitely near points have a large history, both classical under the Zariski school, and modern,
with blowing-up and sheaf theoretic tools.
However, the objective here is to introduce a seemingly bona fide homogeneous ideal {\em on the nose}
in the polynomial ring $k[x,y,z]$ -- as is the case of the ordinary fat ideal --
that contains and is closely related
to the base ideal.
Besides, it comes along with another homogeneous ideal  defined in terms of divisorial valuations; as it turns out,
the latter is a nice carrier to questions about the integral
closure $\bar{I}$ of the base ideal $I$.
Including the usual fat ideal (associated to the proper base points of the map), one finds the base ideal as
a successive subideal of three more ideals.
It is possible moreover to bring into the picture both the saturation of $I$ and its integral closure $\bar{I}$.
In the case of a Cremona map, we end up with a sequence of inclusions with $I$ and the fat ideal as extremes
and yields equalities throughout of the corresponding linear systems spanned in the degree of the rational map.
We intend to explore better the inclusions between these ideals in a future work.

An additional aspect of interest, as a consequence of the present ideal theoretic and homological steps,
is a short proof that a Cremona map of degree $\leq 4$
is saturated (Theorem~\ref{from_lower_bound} (i)), a result that does not seem to have been explicitly given before.

\smallskip

The second section studies the specifics of Cremona maps. First is the case of de Jonqui\`eres maps, whose crucial
role is well-known.
De Jonqui\`eres maps generate a subgroup of the entire group of Cremona transformations with well-known group structure.
 Recently, there has been some intensive activity around the dynamic of these maps.
 The latter catches the asymptotic character of the iterates of a Cremona map.
 Our purpose here is to convey the ideal theoretic properties of a given  map -- a {\em static} view instead.

 We introduce basic ideal theoretic properties of de Jonqui\`eres map
(Corollary~\ref{jonquieres}), including a neat description of the syzygies and the Rees ideal (defining the blowup)
of the base ideal, in addition to an analysis of the behavior of the powers of the latter (Theorem~\ref{Rees_of_Jonq}).
Such a complete picture of the ideal theoretic side of these maps does not seem to be found in the previous literature.
In any case, it is a first step toward a classification of homaloidal types.
A main result of this part is Theorem~\ref{tcharactrization} which explains the homological details of
non-saturated base ideals in degrees $5,6,7$.
The method can in principle be used to classify higher degrees as well, but the results become more involved to
describe.

Since $5$ is the first degree where the base ideal of a plane Cremona map is not necessarily saturated,
we carry a homological characterization of the three known proper homaloidal types in this degree.
The result is given in Theorem~\ref{proper_types_degree5}.

At the end of the section we state a homological criterion for a restricted class of rational maps in
degree $4$ to be Cremona (Proposition~\ref{ideal_squared_replaces_telescopic}).
This trades  the Cohen--Macaulayness
hocus-pocus for the geometric background. This is just an example of such theorems and we believe that there are more of the kind.
A couple of questions is stated that may suggest further work on the subject along the present lines.

\medskip

For convenience we next give a short recap of the terminology used in classical plane Cremona map theory (see \cite{alberich}
for most of the notions to follow).

Let $F\colon
\pp^2\dasharrow \pp^2$ be a rational map. Knowingly, $F$ is defined by three forms of the same degree $d\geq 1$.
Its  {\em base ideal} is the ideal $I\subset R=k[x,y,z]$ generated by these forms.
As is of wide acceptance, the $k$-vector subspace of $R_d$ spanned by these forms is often called a {\em linear system} on $\pp^2$
(whereas it is in fact a space of global sections of a suitably defined linear system or locally free sheaf of rank
one on $\pp^2$).
Clearly, $I$ is generated by the elements of this linear system which allows us to abuse going back and forth between the
two, often in an imprecise style.

If these forms have no nontrivial common factor -- i.e., if the linear system has no fixed part -- then $F$ is said
to have {\em degree} $d$.
In this case, $I$ has codimension $\geq 2$ and hence, the radical of $I$ defines a reduced finite set $V(I)$ of
points. Clearly, the codimension of $I$ is in this case exactly $2$ if and only if $I$ is not $R_+$-primary
(equivalently,  the set of base points is nonempty).
The points of $V(I)\subset \pp^2$ are called {\em proper base points} of $F$ -- this terminology is suggestive
of the existence of other ``improper'' base points.
Thus, if $I$ has codimension $2$ then it has a primary decomposition whose minimal primary components are associated
to the minimal primes of $R/I$ defining the proper base points
$V(I)=\{p_1,\ldots,p_r\}$.
Since $k$ is algebraically closed, every one of these primes is generated by two independent linear forms.

We next recall the notion of (effective) multiplicity.
Namely, given a variety $X$, a smooth point $p\in X$ and a hypersurface (divisor) $D$ then the {\em multiplicity} $e_p(D)$ of $D$ at $p$
is the order of vanishing of a local equation of $D$ at $p\,$; algebraically, if $f$ is a local equation of $D$ at $p$,
then $e_p(D)=\min\{r\geq 0\,|\, f^r\in \mathfrak{m}\}$, where $ \mathfrak{m}$ is the maximal ideal of the local ring
of $X$ at $p$.
For our purpose, $X$ will always be a smooth projective surface; more particularly, $X$ is either $\pp^2$ or the resulting
surface of repeatedly blowing up a reduced finite set of points on the latter.

Going back to our rational map $F\colon
\pp^2\dasharrow \pp^2$, for each
proper base point $p_j$, $j=1,\ldots,r$, one introduces into the picture a {\em virtual multiplicity}
$\mu_{p_j}=\mu_{p_j}(F):=\min\{ e_{p_j}(f)\,|\,f\in I_d\}$.
Note that the subset of the linear system whose elements have  at $p_j$ effective multiplicity equal
to $\mu_{p_j}$ forms an open set $U_{p_j}$ (in the set of parameters).

The complication in the theory is that, besides proper base points, $F$ has built in
other base points, called infinitely near base points.
To define these points, fix a proper base point $p$ of $F$ and let $\mathcal{B}_p$ denote
the surface obtained by blowing-up $p$ on $\pp^2$. Let $E_p\subset \mathcal{B}_p$ stand for
the corresponding exceptional divisor.
A point $q\in E_p$ is a {\em base point infinitely near} to $p$ if it belongs to the proper transform of every divisor in
the open set $U_p$ as introduced above.
We then define the {\em virtual multiplicity} of $F$ at $q$ to be the integer
$\mu_q=\mu_q(F):=\{e_q(\tilde{D})\,|\, D\in U_p\}$, where $\, \tilde\,$ denotes proper transform.
This procedure is to be repeated successively (see \cite[Definition 1.1.49]{alberich} for
the details).

The total set $K$ of points, proper and infinitely near ones, is called the (the set of) {\em base points} of $F$.
The complete set $\boldsymbol{\mu}=\{\mu_1,\ldots,\mu_r\}$ of  multiplicities corresponding to the
base points, preceded by the degree $d$ of $F$, is the
{\em characteristic} (\cite{alberich}) of $F$ -- usually denoted $(d\,; \mu_1,\ldots,\mu_r)$
or $(d\,; \mu_1^{m_1},\ldots,\mu_n^{m_n})$ to account for the repetitions.
The resulting {\em cluster} $\mathcal{K}=(K, \boldsymbol{\mu})$ is the (weighted) {\em cluster of base points} of $F$
or of the corresponding linear system -- a notion that plays an important
 role in the classical Cremona theory; we refer
to \cite{Casas} which contains a detailed study of this notion and its role in plane curve theory.
Thus, any rational map $F:\pp^2\dasharrow \pp^2$ with codimension $2$ base ideal, and in particular, a Cremona map
with  codimension $2$ base ideal, carries a weighted cluster
$\mathcal{K}=(K,\boldmath{\mu})$, consisting of the
set of its base points suitably ordered, along with the set of their corresponding virtual multiplicities
(weights) as explained above.

We also recall the classical  {\em equations of condition} for a plane Cremona map  of
degree $d$ (see \cite[2.5]{alberich}):
\begin{equation}\label{eqs_condition}
\sum_{p}\mu_p=3d-3,\; \sum_{p}\mu_p^2=d^2-1,
\end{equation}
where $p$ runs through the set of (proper and infinitely near) base points of the corresponding
linear system with respective  multiplicities $\mu_p$.
An abstract configuration $(d\,;\mu_1,\ldots,\mu_r)$ satisfying the equations of condition
is called a {\em homaloidal type}, and is denoted in the same fashion as the characteristic of
a Cremona map.

A Cremona map whose base points are proper is called {\em simple}.

A homaloidal type is called {\em proper} if there exists a plane Cremona map whose characteristic
coincides with it. There is an important practical  tool to test whether a given homaloidal
type is proper - it is called {\em Hudson test} (\cite[Corollary 5.3.2]{alberich}).

\medskip

{\bf Acknowledgements.} The first author is grateful to M. Alberich-Carami\~nana, C. Ara\'ujo,
M. Chardin and I. Pan for helpful discussions.
The second author thanks C. Ciliberto and F. Russo for an illuminating conversation at the early
stages of the work.

\section{Related ideal theory}\label{first}

\subsection{Homological results}

In this part we develop some generalities on resolutions
of ideals generated by three forms.

The following basic statement about graded minimal resolutions does not seem to have been noted before
 in this particular form.
Its only use will be in the standard graded case, but because of its possible independent interest
it may be convenient to state it in the larger realm of positively graded Noetherian $^*$local rings
(see \cite[Definition 1.5.13]{BHbook}).

\begin{Lemma}\label{lgeneral}Let $(R,\fm)$ be a positively graded Noetherian $^*$local ring and let
$I\subset \fm$ be a $3$-generated homogeneous ideal  with graded  minimal free resolution
\begin{equation}\label{eqresolution}
0\rightarrow\bigoplus_{m=1}^{r-2}R(-D_m)\xrightarrow{\phi_3}\bigoplus_{i=1}^{r}R(-d_i)
\xrightarrow{\phi_2}\bigoplus_{l=1}^{3}R(-a_l)\xrightarrow{\phi_1}R\rightarrow R/I \rightarrow 0.
\end{equation}
Assume $D_1\geq\cdots\geq D_{r-2}$ and $d_1\geq\cdots\geq d_{r}$.
Then $D_m\geq d_m +1$ for all $1\leq m \leq r-2$.
\end{Lemma}
\demo By the Buchsbaum--Eisenbud criterion \cite[Theorem 20.9]{E}, $\rank(\phi_3)=r-2$ and $\grade(I_{r-2}(\phi_3))\geq 3$.
Since the resolution is graded minimal, $\grade(I_{r-2}(\phi_3))= 3$.
Therefore, the Eagon--Northcott is a graded minimal resolution of $I_{r-2}(\phi_3)$,
which implies in particular that the maximal minors of $\phi_3$ form a minimal set of generators of the ideal generated by these minors.
In particular, every maximal minor is nonzero.

Now set $\phi_3=(m_{ij})$, with $m_{ij}\in\fm,\,\forall i,j$.  We claim that $\deg(m_{ii})\geq 0$ for all $1\leq i \leq r-2$.
Indeed, if $D_i\leq d_i$ for some $i$ then for all $p \leq i$ and $q \geq i$ we would have
$\Deg(m_{pq})=D_q-d_p\leq D_i-d_i\leq 0$ and since $m_{pq}\in\fm$, necessarily $m_{pq}=0$ for all $p \leq i$ and $q \geq i$.
It then follows that the upper $(r-2)\times(r-2)$ submatrix  of $\phi_3$ is of the following shape:
$$
\left(
\begin{array}{cccccl}
*&&..&0     &\cdots&0\\
 &*&..&0&\cdots&\vdots\\
 &&     &     0&\cdots&0\\
  &&&&*&\\
 &&    &     &\cdots&*
\end{array}
\right)
$$
This determinant is obviously null, thus  contradicting the above assertion. Therefore $\Deg(D_i)-\Deg(d_i)=\Deg(m_{ii})>0$ for all $1\leq i \leq r-2$.
\qed

\bigskip


Our focus in this work is on the case where $(R,\fm)=(k[x_0,x_1,x_2], (x_0,x_1,x_2))$
in which $k$ is a field and $R$ is assigned the standard grading.
We set $M\,\Check {}$ for the graded Matlis dual of an $R$-module into $k=R/\fm$.

For a graded $R$-module $M$, we will denote
$$\beg (M):=\inf \{ \mu \ \vert \ M_\mu \not= 0\},$$
with the convention that $\beg (0)=+\infty$, and
$$\END (M):= \sup \{ \mu \ \vert \ M_\mu \not= 0 \},$$
with the convention that $\END(0)=-\infty$.

As a matter of further notation, we set $I^{\rm sat}:=I:\fm^{\infty}$
and denote by $I^{\rm un}$ the unmixed part of the primary decomposition
of the ideal $I$.
Finally,  $\omega_{R/I}$ will denote the graded canonical module
of $R/I$.

The duality piece in the next result is a special case of  \cite[Lemma 5.8]{Ch}, but we
will give a proof for the reader's convenience and later reference.

\begin{Theorem}\label{ldual}
Let $R=k[x_0,x_1,x_2]$ and let $I\subset R$ be an  ideal of height $2$ generated by $3$ linearly
independent forms of degree $d\geq 1$. Then
 \begin{itemize}
\item[{\rm (i)}] $(I^{\rm sat}/I)\Check {}\simeq (I^{\rm sat}/I) (3d-3)\,${\rm ;} in particular,
 if $I$ is not saturated then
$$\END(I^{\rm sat}/I)+\beg(I^{\rm sat}/I)=3d-3.$$
\item[{\rm (ii)}] $\END(H^1_{\fm}(R/I))+1=-\beg(\omega_{R/I})+1\leq 2d-3$.
 \end{itemize}
\end{Theorem}
\demo
Note at the outset that $I$ is a strict almost complete intersection generated in degree $d\geq 2$.

Consider the graded Koszul complex $K_{\bullet}$ generated by a minimal generating set of $I$ and set $H_i:=H_i(K_{\bullet})$ for its $i$th homology. Consider the complexes tensor product $K_{\bullet}\otimes C^{\bullet}_{\fm}$, where $C^{\bullet}_{\fm}$ is the \v {C}ech complex on $\fm$.
Consider the spectral sequence $\E$ associated to the total complex of the double complex $K_{\bullet}\otimes C^{\bullet}_{\fm}$. (For details on this terminology, we refer to \cite[Sections A3.13.4 and A3.13.5]{E}.)
Note that $H_{q-p}($Tot$(E^{\bullet,\bullet}))=\sideset{^{\infty}}{_{\text\footnotesize{hor}}^{-p,-q}}{\E}$.
Putting this complex  in the third quadrant, the terms of the second spectral sequence are
$$
\sideset{^2}{_{\text\footnotesize{\rm ver}}^{-p,-q}}{\E}= \left\lbrace
           \begin{array}{c l}
               H^q_{\fm}(R/I)  & \text{if $p=0$ and $q=0,1$},\\
               H^1_{\fm}(H_1)  & \text{if $p=1$ and $q=1$},\\
               0                & \text{otherwise};
           \end{array}
         \right.
$$
$$
\sideset{^2}{_{\text\footnotesize{\rm hor}}^{-p,-q}}{\E}= \left\lbrace
           \begin{array}{c l}
               (H_0(3d-3))\Check {}  & \text{if $p=3$ and $q=3$},\\
               (H_1(3d-3))\Check {}  & \text{if $p=2$ and $q=3$},\\
               0                & \text{otherwise}.
           \end{array}
         \right.
$$

 Therefore, by  the convergence of the spectral sequence, the homology of the total complex is filtered by elements on the  diagonal of  $\sideset{^{\infty}}{_{\text\footnotesize{\rm ver}}^{\bullet,\bullet}}{\E}$ which contains
 $\sideset{^{\infty}}{_{\text\footnotesize{\rm ver}}^{-p,-q}}{\E}$. The only non-trivial filtration thus
  obtained is for $H_0($Tot$(E^{\bullet,\bullet}))=(H_0(3d-3))\Check {}\,$, which is the following short exact sequence:
\begin{equation}\label{homology_sequence}
0 \rightarrow H^1_{\fm}(H_1)\lar (H_0(3d-3))\Check {} \lar H^0_{\fm}(R/I) \lar 0.
\end{equation}

Since $I$ is an almost complete intersection, one has $H_1\simeq \omega_{R/I}\simeq \Ext ^2(R/I,R)(-3d)$. Moreover,  quite generally, $\omega_{R/I}=\omega_{R/I^{\rm un}}$. Therefore, by graded duality  (\cite{BHbook})
\begin{center} $H^1_{\fm}(H_1) \simeq (\Ext^2(\Ext^2(R/I^{\rm un},R)(-3d),R(-3)))\Check {}\simeq (R/I^{\rm un}(3d-3))\Check {}$,
\end{center}
where  the rightmost isomorphism is due to the Cohen-Macaulayness of $R/I^{\rm un}$.

Then the exact sequence (\ref{homology_sequence}) becomes
\begin{equation}\label{duality_sequence}
0 \rightarrow (R/I^{\rm un}(3d-3))\Check {}\lar (R/I(3d-3))\Check {} \lar H^0_{\fm}(R/I) \lar 0.
\end{equation}

Thus $H^0_{\fm}(R/I)\simeq (I^{\rm un}/I)\Check {}\,(3-3d)$. Since $H^0_{\fm}(R/I)= I^{\rm sat}/I$ and  $I^{\rm sat}=I^{\rm un}$
in the present context, we are done for (i).

To prove (ii), note that $R/I^{\rm sat}$ is a CM ring of dimension $1$ and that  $I^{\rm sat}/I$ is a module of
finite length; thus $H^1_{\fm}(R/I)=H^1_{\fm}(R/I^{\rm sat})$. Then
\begin{eqnarray*}
\END(H^1_{\fm}(R/I))+1&=&-\beg(\omega_{R/I^{\rm un}})+1 = -\beg(\omega_{R/I})+1\\
&=& -\beg(\Ext ^2(R/I,R)(3d-3))+1.
\end{eqnarray*}
 On the other hand, we may assume that $k$ is infinite, hence
 $$\Ext ^2(R/I,R)(3d-3)=(\alpha:I/\alpha)(2d-3),$$
 where $\alpha$ is a maximal regular sequence of  $d$-forms in $I$ .
 Since $I$ is not complete intersection,  $\beg((\alpha:I/\alpha))\geq 1$. Collecting the information, we arrive at the
statement.
\qed

\begin{Proposition}\label{gen_degs_of_H0} Let $R=k[x_0,x_1,x_2]$ and let $I\subset R$ be an ideal of codimension $2$
generated by $3$ linearly independent forms
of degree $d\geq 1$ with minimal graded free resolution
\begin{equation}\label{eqdresolution}
0\rightarrow\bigoplus_{i=1}^{r-2}R(-D_i)\longrightarrow \bigoplus_{i=1}^{r}R(-d_i)\longrightarrow
R^{3}(-d)\rightarrow R\rightarrow R/I \rightarrow 0 \quad (r\geq 3).
\end{equation}
Then:
\begin{itemize}
\item[{\rm (i)}] The minimal free resolution of $I^{\rm sat}/I$ as an $R$-module has the form
{\small
\begin{equation}\label{eqdresolutionIsatI}
0\rar\bigoplus_{i=1}^{r-2}R(-D_i)\rar\bigoplus_{i=1}^{r}R(-d_i)\rar \bigoplus_{i=1}^{r}R(d_i-3d)\rar \bigoplus_{i=1}^{r-2}R(D_i-3d)
\rightarrow I^{\rm sat}/I \rightarrow 0
\end{equation}
}
where the leftmost map is the same as that of {\rm (\ref{eqdresolution})}.
\item[{\rm (ii)}] If in addition $I^{\rm sat}_d=I_d$ and $I^{\rm sat}_i=0$ for $i< d$, then the resolution
of $I^{\rm sat}$ is
\begin{equation}\label{eqdresolutionIsat}
0\longrightarrow \bigoplus_{i=1}^{r}R(-(3d-d_i))\longrightarrow R^{3}(-d)\bigoplus_{i=1}^{r-2}R(-(3d-D_i))
\rightarrow  I^{\rm sat} \rightarrow 0.
\end{equation}
\end{itemize}
\end{Proposition}
\demo (i)
Applying $\Hom_R(-,R(-3))$ to (\ref{eqdresolution}) yields a
minimal free presentation
$$\bigoplus_{i=1}^{r}R(d_i-3)\longrightarrow \bigoplus_{i=1}^{r-2}R(D_i-3)\longrightarrow \Ext^3(R/I,R(-3)) \rightarrow 0.$$
On the other hand, by graded duality and by the exact sequence (\ref{duality_sequence}) one has
$$\Ext^3(R/I,R(-3))\simeq H^0_{\fm}(R/I)\Check {}=(I^{\rm sat}/I)(3d-3).$$
Therefore, we have a free presentation
\begin{equation}\label{presentation_of_quotient}
\bigoplus_{i=1}^{r}R(d_i-3)\longrightarrow \bigoplus_{i=1}^{r-2}R(D_i-3)\longrightarrow (I^{\rm sat}/I)(3d-3) \rightarrow 0.
\end{equation}
Shifting by $-3d+3$ yields a free resolution of the form
\begin{equation}\label{eqdresolutionMef}
0\rightarrow\bigoplus_{i=1}^{l-2}R(-a_i)\longrightarrow\bigoplus_{i=1}^{l}R(-b_i)\longrightarrow
\bigoplus_{i=1}^{r}R(d_i-3d)\longrightarrow \bigoplus_{i=1}^{r-2}R(D_i-3d)\rightarrow I^{\rm sat}/I \rightarrow 0,
\end{equation}
with suitable integers $l, a_i, b_i$.
Applying $\Hom_R(-,R(-3))$ to the latter yields a free complex resolving the third homology $\Ext^3(I^{\rm sat}/I, R(-3))$.
But, again by graded duality plus the fact that $H^0_{\fm}(I^{\rm sat}/I)=I^{\rm sat}/I$ by definition of $I^{\rm sat}$, and
using Theorem~\ref{ldual} (i), we get
$$\Ext^3(I^{\rm sat}/I, R(-3))\simeq (I^{\rm sat}/I)(3d-3).$$
Shifting by $-3d+3$ once more yields a free resolution of the form
\begin{equation}\label{eqdresolutionM}
0\rightarrow\bigoplus_{i=1}^{r-2}R(-D_i)\longrightarrow\bigoplus_{i=1}^{r}R(-d_i)\longrightarrow
\bigoplus_{i=l}^{l}R(b_i-3d)\longrightarrow \bigoplus_{i=1}^{l-2}R(a_i-3d)\rightarrow I^{\rm sat}/I \rightarrow 0.
\end{equation}
Comparing (\ref{eqdresolutionMef}) and (\ref{eqdresolutionM}), the uniqueness of the minimal free resolution yields the assertion.

\smallskip

(ii) Since $R/I^{\rm sat}$ is Cohen--Macaulay, the hypothesis of the item implies a graded minimal resolution of the form
$$0\longrightarrow \bigoplus_{i=1}^{s+2}R(\alpha_i)\longrightarrow \left(\bigoplus_{i=1}^{s}R(-\beta_i)\right)\oplus R^{3}(-d)
 \rightarrow  I^{\rm sat} \rightarrow 0.$$
 In order to determine the shifts $\alpha_i, \beta_i$, we consider a map of complexes lifting the inclusion $I\subset I^{\rm sat}$
 {\small
 $$\begin{array}{cccccccclccc}
0&\rightarrow  &0 &\rightarrow &\bigoplus_{i=1}^{s+2}R(\alpha_i)&\rightarrow &
\left(\bigoplus_{i=1}^{s}R(-\beta_i)\right)\oplus R^{3}(-d))&\rightarrow & I^{\rm sat}&\rightarrow&0& \\[5pt]
 &    & \uparrow  & &   \uparrow   &   &\stackrel{\iota}{}\;\uparrow & &\uparrow &  & &\\[5pt]
0&\rightarrow &\bigoplus_{i=1}^{r-2}R(-D_i)&\rightarrow &\bigoplus_{i=1}^{r}R(-d_i)  &\rightarrow &  R^{3}(-d) &\rightarrow &I &\rightarrow&0&
\end{array}$$
}
with $\iota$ the inclusion in the ``second'' coordinate.
Now, the mapping cone of this map
{\small
$$ 0 \rar \bigoplus_{i=1}^{r-2}R(-D_i)\lar \bigoplus_{i=1}^{r}R(-d_i)\lar \bigoplus_{i=1}^{s+2}R(\alpha_i)\bigoplus  R^{3}(-d)
\lar (\bigoplus_{i=1}^{s}R(-\beta_i))\oplus R^{3}(-d)$$
}
is a free resolution of $I^{\rm sat}/I$.
Canceling the non-minimal part coming from the summand  $R^{3}(-d)$, yields a minimal free complex.
Comparing with (\ref{eqdresolutionIsatI}) yields $s=r, \alpha_i=d_i-3d, \beta_i=3d-D_i$.
\qed

\subsection{A critical lower bound for $I^{\rm sat}/I$}

If $I$ is a homogeneous ideal minimally generated in some degree $d\geq 1$, but non-saturated,
 typically $I^{\rm sat}/I$ will have minimal generators in degrees $\leq d$.
In this part, we derive substantial consequences by requiring that the initial
degree of $I^{\rm sat}/I$ be $d+1$.
For later use, we remark that this critical bound holds true when $I$ is the
base ideal of a Cremona map (\cite[Proposition 1.2]{PanRusso}).

First is an upper bound for the regularity in terms of $d$.

\begin{Corollary}\label{strong_bound} Let $R=k[x_0,x_1,x_2]$ and let $I\subset R$ be an ideal of
codimension $2$ generated by $3$ linearly independent forms of degree $d\geq 1$.
If  $\,\beg(I^{\rm sat}/I)\geq d+1$ then $\Reg(R/I)\leq  2d-3$.
\end{Corollary}
\demo
By definition, one has
$$\Reg(R/I)=\MAX\{\END(H^0_{\fm}(R/I)), \END(H^1_{\fm}(R/I))+1\}.$$
By Theorem~\ref{ldual}(ii),  $\END(H^1_{\fm}(R/I))+1\leq 2d-3$.
Thus, we are done if $I$ is saturated because $\END(H^0_{\fm}(R/I))=-\infty$.

If $I$ is not saturated, then
$$\END(H^0_{\fm}(R/I))=3d-3-\beg(I^{\rm sat}/I)\leq 3d-3-(d+1)=2d-4,$$
by Theorem~\ref{ldual}(i) and the hypothesis.
Hence we are done in this case too.
\qed

\begin{Theorem}\label{from_lower_bound} Let $R=k[x_0,x_1,x_2]$ and let $I\subset R$ be an ideal
of codimension $2$ generated by $3$ linearly independent forms
of degree $d\geq 1$ with minimal graded free resolution
\begin{equation}\label{eqdresolution2}
0\rightarrow\bigoplus_{m=1}^{r-2}R(-D_m)\longrightarrow \bigoplus_{i=1}^{r}R(-d_i)
\longrightarrow R^{3}(-d)\rightarrow R\rightarrow R/I \rightarrow 0 \quad (r\geq 3).
\end{equation}
If $\,\beg(I^{\rm sat}/I)\geq d+1$ then:
\begin{itemize}
\item[\rm (i)]{$d\geq 5$}
\item[\rm (ii)]{$\Reg(R/I)=3d-3-\beg(I^{\rm sat}/I)=D_1-3\leq 2d-4$}
\item[\rm (iii)]{$r \leq d-2$.}
\end{itemize}
\end{Theorem}
\demo
(i)$\,$ This follows from the equation  $\END(I^{\rm sat}/I)+\beg(I^{\rm sat}/I)=3d-3$
in Theorem~\ref{ldual}(i) and from the assumption $\beg(I^{\rm sat}/I)\geq d+1$.
 Indeed, one has
 $$2d+2\leq 2\beg(I^{\rm sat}/I)\leq \END(I^{\rm sat}/I)+\beg(I^{\rm sat}/I)=3d-3,$$
 hence $d\geq 5$.

 \smallskip

 (ii)$\,$ Since $\Reg(R/I)=\MAX\{d-1,d_1-2,D_1-3\}$, Lemma~\ref{lgeneral} gives $\Reg(R/I)=D_1-3$.
 By Proposition~\ref{gen_degs_of_H0} and the hypothesis, $3d-D_1=\beg(I^{\rm sat}/I)\geq d+1$.
 Assembling yields $\Reg(R/I)=D_1-3=3d-3-\beg(I^{\rm sat}/I)\leq 2d-4$.

 \smallskip

 (iii)$\,$ From the resolution of $R/I$ its Hilbert series is
 \begin{equation}\label{hseries}
\frac{1-3t^d+\sum_{i=1}^r t^{d_i}- \sum_{m=1}^{r-2} t^{D_m}}{(1-t)^3},
\end{equation}
with a pole of order $2$ at $t=1$ since $\dim R/I=1$.
Taking $t$-derivatives of the numerator of (\ref{hseries}) evaluated at $t=1$
(see \cite[4.1.14]{BHbook}), one obtains
the following relation
 \begin{equation}\label{firstderivative}
d_r+d_{r-1}=\sum_{m=1}^{r-2}( D_m -d_m)+3d.
\end{equation}
 Now, Lemma \ref{lgeneral} implies that $\sum_{m=1}^{r-2}( D_m -d_m)+3d\geq r-2+3d$.
 On the other hand, part (ii) yields $\Reg(R/I)\leq 2d-4$ which implies $d_1\leq 2d-2$.
Therefore $d_r+d_{r-1}\leq 4d-4$.
Assembling the inequalities, we get $r-2+3d\leq 4d-4$, hence $r\leq d-2$ as required.
\qed

\medskip

As a consequence of the above results, one can classify the virtual resolutions
along with the corresponding twists in case of low values of $d$.
In the subsequent sections we will deal with the question whether such virtual
resolutions are in fact realized by the base ideal of a plane Cremona map -- note that realisability
is a question only when the resolution has length at least $3$.
We emphasize that in the case where $I$ is saturated the virtual resolutions are easily computed by
 means of (\ref{firstderivative}) and Corollary~\ref{strong_bound}.

\begin{Proposition}\label{P567}
With the same notation and hypotheses as in {\rm Theorem~\ref{from_lower_bound}}, one has:
\begin{itemize}
\item[\rm (i)]{If $d=5$, the minimal free resolution of $R/I$ is \\ $0\rightarrow R(-9)\rightarrow R^3(-8)\rightarrow R^3(-5)\rightarrow R$}
\item[\rm (ii)]{If $d=6$, the minimal free resolution of $R/I$ is one of the following:
\\$0\rightarrow R^2(-11)\rightarrow R^4(-10)\rightarrow R^3(-6)\rightarrow R$, or
\\$0\rightarrow R(-11)\rightarrow R^2(-10)\bigoplus R(-9) \rightarrow R^3(-6)\rightarrow R$.}
\item[\rm (iii)]{If $d=7$, the minimal free resolution of $R/I$ is one of the following:
\\$0\rightarrow R^3(-13)\rightarrow R^5(-12)\rightarrow R^3(-7)\rightarrow R$,
\\$0\rightarrow R^2(-13)\rightarrow R^3(-12)\bigoplus R(-11) \rightarrow R^3(-7)\rightarrow R$,
\\$0\rightarrow R(-13)\rightarrow R(-12)\bigoplus R^2(-11) \rightarrow R^3(-7)\rightarrow R$,
\\$0\rightarrow R(-12)\rightarrow R^3(-11) \rightarrow R^3(-7)\rightarrow R$, or
\\$0\rightarrow R(-13)\rightarrow R^2(-12)\bigoplus R(-10) \rightarrow R^3(-7)\rightarrow R$
.}
\end{itemize}
\end{Proposition}
\demo
First recall the following strand of inequalities from the proof of
Theorem~\ref{from_lower_bound}, (iii):
\begin{equation}\label{ineq}
4d-4\geq d_r+d_{r-1}=\sum_{m=1}^{r-2}( D_m -d_m)+3d\geq r-2+3d\geq 1+3d.
\end{equation}

(i) $\,$ If $d=5$, then $r=3$ by Theorem~\ref{from_lower_bound}, (iii) and the inequalities in
(\ref{ineq}) are all equalities, whence $d_3=d_2=d_1=2d-2=8$ and $D_1-d_1=1$.

\smallskip

(ii) $\,$ If $d=6$ then $r=3$ or $r=4$ by Theorem~\ref{from_lower_bound} (iii).
In the case where $r=4$, the inequalities in (\ref{ineq}) are all equalities except the
rightmost one and there is one single solution. The shifts turn out to be as stated.
If $r=3$ then (\ref{ineq}) yields two virtual solutions, with $D_1=11$ or $D_1=12$.
However, by Theorem~\ref{from_lower_bound} (ii), only $D_1=11$ is possible and the
shifts are as stated.

\smallskip

(ii) $\,$ If $d=7$ then $3\leq r\leq 5$ by Theorem~\ref{from_lower_bound} (iii).
The case $r=5$ yields again equalities throughout
(\ref{ineq}) except for the rightmost inequality and there is one single solution, as stated.

For $r=4$, from  (\ref{ineq}) we deduce that $d_3\geq r+d+1=12=2d-2$, hence $d_3=d_2=d_1=2d-2=12$
as this is the highest possible value.
As a consequence, both $D_1$ and $D_2$ attain the upper bound $2d-1=13$.
Thence the equality $d_4+d_3=\sum_{m=1}^{r-2}( D_m -d_m)+3d=23$ gives $d_4=11$, as stated.

For $r=3$ the argument is slightly more involved, but again (\ref{ineq}) yields the result: the third and fourth
possibilities on the list stem from the equality $d_2=d+r+1=11$ which gives $d_3=11$ and $d_1=12$ or
$d_1=11$ and accordingly $D_1=d_1+1=13, 12$.
The last case listed is when $d_2=12$, the largest possible value. Accordingly, from (\ref{ineq})
comes $d_2=d_1=D_1-1=12$, hence $d_3=10$.
\qed

\medskip

For an  ideal $I$ in a standard graded or local ring $(R,\fm)$ over a field one defines the {\em saturation exponent}  of $I$ as
${\rm st}(I):=\MIN\{s\in \mathbb{N}\,|\, I^{\rm sat}=I:\fm^s\}$.
This number has been variously called {\em saturation index} or {\em satiety}.

\begin{Corollary}\label{psat}
With the same notation and hypotheses as in {\rm Theorem~\ref{from_lower_bound}}, one has:
\begin{itemize}
\item[{\rm (i)}] {\rm st}$(I)\leq 2\Reg(R/I)-3d+4\leq d-4$
\item[{\rm (ii)}] If the module $I^{\rm sat}/I$ is minimally generated by elements of a fixed degree
then {\rm st}$(I)= 2\Reg(R/I)-3d+4$. In particular, this equality holds for $d\leq 7$.
\end{itemize}
\end{Corollary}
\demo
(i)$\,$ Let $i_0=\beg(I^{\rm sat}/I)$ and $e=\END(I^{\rm sat}/I)$.
Then
$$\fm^{e-i_0+1}(I^{\rm sat}/I)_j\subset (I^{\rm sat}/I)_{e+j-i_0+1}=\{0\},$$
for all $j\geq 0$, because $j\geq i_0$ as soon as $(I^{\rm sat}/I)_j\neq \{0\}$.
Therefore st$(I)\leq e-i_0+1$. On the other hand Theorem~\ref{ldual} (iii) and
Theorem~\ref{from_lower_bound} (ii)
yield $e-i_0+1=2\Reg(R/I)-3d+4\leq d-4$.

(ii)$\,$ By Proposition~\ref{gen_degs_of_H0},    $I^{\rm sat}/I$ is generated by
elements in $(I^{\rm sat}/I)_{i_0}$. Therefore it obtains  $0\neq (I^{\rm sat}/I)_e\subseteq \fm^{e-i_0}(I^{\rm sat}/I)_{i_0}$
thus yielding the required equality.

Finally, for $d\leq 7$, Proposition~\ref{P567} together with Proposition~\ref{gen_degs_of_H0}(i) shows that $I^{\rm sat}/I$ is generated in fixed degree.
\qed

\subsection{Analogues of the fat ideal}\label{analogues}

In this portion we wish to get hold of other homogeneous ideals closely related to the base ideal
of a rational map $F:\pp^2\dasharrow\pp^2$ without fixed part.
For the sake of clarity, we briefly go back to some of the notions given in the Introduction.

Letting $\mathcal{K}=(K,\boldsymbol\mu)$
denote the weighted cluster of the linear system defining $F$, one can consider the {\em $\boldsymbol\mu$-fat ideal}
$$\mathfrak{F}_{\mathcal{K}}:=I(p_1)^{\mu_1}\cap\cdots\cap I(p_r)^{\mu_r},$$
corresponding to the proper points of $K$, where  $I(p_j)$ stands for the defining prime ideal of the proper point $p_j$.
Thus, for a plane rational map $F$ whose base points are all proper, the $\boldsymbol\mu$-fat ideal is a rough ideal theoretic
saturated approximation to the base ideal of $F$.

One advantage of the fat ideal is that the degree of the corresponding scheme is automatic:
\begin{equation}\label{multplicity_of_fat}
e(R/\mathfrak{F}_{\mathcal{K}})=\sum_{j=1}^r\, \frac{\mu_j(\mu_j+1)}{2}.
\end{equation}
Now, if $k$ is algebraically closed (or if we consider only rational points) then
$I\subset \mathfrak{F}_{\mathcal{K}}$. This is because saying that a form $f\in R$ has multiplicity
at least $\mu_j$ locally at the prime $I(p_j)\subset R$ is equivalent to asserting that
$f\in I(p_j)^{(\mu_j)}=I(p_j)^{\mu_j}$ (the equality holding because $I(p_j)$ is a complete intersection).
This allows for a comparison of the two ideals and shows that $e(R/I)\geq  \sum_{j=1}^r\,
 \frac{\mu_j(\mu_j+1)}{2}$.

On the other hand, since the saturation $I^{\rm sat}=I:(R_+)^{\infty}$ of $I$  coincides with
the unmixed part of $I$, it is the smallest
unmixed ideal (by inclusion) containing $I$.
Since $\mathfrak{F}_{\mathcal{K}}$ is unmixed by definition and has same radical as $I$, one has the valuable setup
$I\subset I^{\rm sat}\subset \mathfrak{F}_{\mathcal{K}}$,  with all three ideals
 sharing the same radical.

Alas, if the map has infinitely near points among its base points, then the $\boldsymbol\mu$-fat ideal is no longer the
 ``tightest'' ideal to look at.
A ``correction'' is available by introducing two new ideals, one of which is based on the classical
theory of (weighted) clusters and the corresponding blowup gadgets.
We now proceed to establish these definitions.

\smallskip

\begin{Definition}\label{passing_virtually}\rm
Let $\mathcal{K}=(K,\boldsymbol\mu)$ denote a weighted cluster.
\begin{itemize}
\item A plane curve $C\subset \pp^2$ {\em passes virtually through}  $\mathcal{K}$ if
the divisor on $\mathcal{B}_K$
\begin{equation}\label{passing_total}
\bar{C}^K-\sum_{p\in K}\, \mu_p\, \bar{E}_p^K
\end{equation}
is effective, where $\mathcal{B}_K$ is the blowup of the set $K$ on  $\pp^2$, $\bar{C}^K$ and $\bar{E}_p^K$ are the
total transform of $C$ and the total $p$-component
of the exceptional divisor $E^K$, respectively (see \cite[Definition 1.1.38]{alberich}.
\item Given $m\geq 1$, let $\ell_{\mathcal{K}}(m)\subset R_m$ consist of all forms $f$ of degree $m$ such the
curve $C=V(f)$ passes virtually through $\mathcal{K}$.
\end{itemize}
\end{Definition}
Set $I_{\mathcal{K}}:=\oplus_{m\geq 0} \ell_{\mathcal{K}}(m)\subset R$.
Then $I_{\mathcal{K}}$ is an ideal (not just a vector subspace) of $R$ - we call it the {\em full ideal of curves
through} $\mathcal{K}$.

The difficulty in handling the algebraic properties of the full ideal of curves is due to the nature of the notion
in (\ref{passing_total}).
One can modify it to requiring, equivalently, that
\begin{equation}\label{passing_proper}
\tilde{C}^K-\sum_{p\in K}\, (e_p(C)-\mu_p)\, \bar{E}_p^K
\end{equation}
be effective, where $\tilde{C}^K$ now denotes the proper transform of $C$ in $\mathcal{B}_K$ (use \cite[lemma 1.1.8]{alberich}).
Since $\tilde{C}^K$ is an effective divisor, the condition that $C$ passes through the cluster $\mathcal{K}$
means that the divisor $\sum_{p\in K}\, (e_p(C)-\mu_p)\, \bar{E}_p^K$ is effective.

If we had $\tilde{E}_p^K$ instead of $\bar{E}_p^K$ in this divisor, then the condition would simply require
the inequalities $e_p(C)\geq \mu_p$ for every $p\in K$.
To still express the actual condition  in terms of inequalities, one resorts to the notion of the {\em proximity
matrix} $\mathbf{P}_{\mathcal{K}}$ of the cluster $\mathcal{K}$, whose entries are exactly the coefficients of each $\tilde{E}_p^K$ in terms
of all $\bar{E}_q^K$, with $q\in K$ (see \cite[Corollary 1.1.27 and Definition 1.1.28]{alberich}).

For reference convenience, we state the final expression in a separate result:
\begin{Lemma}\label{passing_lemma}
A plane curve $C\subset \pp^2$ passes virtually through  $\mathcal{K}$ if and only if
\begin{equation}\label{passing_as_proximity}
\mathbf{P}_{\mathcal{K}}^{-1}\cdot (\mathbf{e}_K(C)-{\boldsymbol\mu}_K)^t\geq 0,
\end{equation}
where $\mathbf{e}_K(C)$ {\rm (}respectively, ${\boldsymbol\mu}_K${\rm )} denotes the vector of effective multiplicities
at the points of $K$ {\rm (}respectively, the vector of virtual multiplicities of $\mathcal{K}${\rm )}, and $^{-1}$, $^t$ denote
matrix inverse and  matrix transpose, respectively.
\end{Lemma}
To make sense, the above inequality relies on the fact that $\det(\mathbf{P}_{\mathcal{K}})\neq 0$ --
actually this matrix is unimodular, so its inverse is an integer matrix; moreover, it can be shown that
the inverse has only nonnegative entries  (\cite[Lemma 1.1.32]{alberich}).

\smallskip

As a preliminary to the subsequent results we can now prove without difficulty:

\begin{Proposition}\label{first_inclusions}
Let $F:\pp^2\dasharrow\pp^2$ be a rational map without fixed part, with cluster of base points
$\mathcal{K}=(K, \boldsymbol\mu)$ and base ideal $I\subset R$.
Then:
\begin{itemize}
\item[{\rm (a)}] $I\subset I_{\mathcal{K}}$.
\item[{\rm (b)}] $I_{\mathcal{K}}\subset \mathfrak{F}_{\mathcal{K}}:=\bigcap_j \wp_j^{\mu_j}$ {\rm (}``fat'' ideal{\rm )},
where $\wp_j\subset R$ is the homogeneous prime ideal of the proper point $p_j\in K$.
If, moreover, $K$ consists of proper points then $I_{\mathcal{K}}=\mathfrak{F}_{\mathcal{K}}$.
\end{itemize}
\end{Proposition}
\demo
(a) Let $d$ be the degree of $F$. By the very nature of the definition of the associated weighted cluster $\mathcal{K}$
of $F$, any curve of the linear system in degree $d$ defining $F$ passes virtually through $\mathcal{K}$ -- note that
the condition in (\ref{passing_as_proximity}) is here verified ``on the nose" through the inequalities $e_p(C)\geq \mu_p$
for every $p\in K$.
Therefore $I_d\subset {(I_{\mathcal{K}})}_d$. It follows that $I_m\subset {(I_{\mathcal{K}})}_m$ for every degree $m\geq d$.
If $m <d$ then $I_m=\{0\}$, hence $I_m\subset {(I_{\mathcal{K}})}_m$ holds trivially.
Thus, $I\subset I_{\mathcal{K}}$.

(b)
$\mathbf{P}_{\mathcal{K}}$ is a block-diagonal matrix and it can be arranged so that the submatrix
corresponding to the proper points is the identity matrix; then the inverse $\mathbf{P}_{\mathcal{K}}^{-1}$
will also contain an identity block that multiplies the proper part of the vector $(\mathbf{e}_K(C)-{\boldsymbol\mu}_K)^t$.
On the other hand, it follows from the definitions that
$$\mathfrak{F}_{\mathcal{K}}=\{f\in R: e_p(f)\geq \mu_p,\; \forall \; {\rm proper}\, p\in K\}.$$
Therefore, the inclusion $I_{\mathcal{K}}\subset \mathfrak{F}_{\mathcal{K}}$ follows immediately.

The second statement of this item is a direct consequence of the part just proved.
\qed

\medskip

To proceed, we need the following order of ideas.

Recall that given a Noetherian domain $R$ with field of fractions $Q$, a valuation ring $(V,\fm_V)$ of $Q$
containing $R$ is called a {\em divisorial valuation ring} on $R$ if the transcendence degree of $V/\fm_V$ over
$R_{\fm_V\cap R}/(\fm_V\cap R)_{\fm_V\cap R}$ is $\hht (\fm_V\cap R)-1$.
The corresponding valuation of $V$ will be called a divisorial valuation relative to $R$.
Let D$(R)$ denote the set of the divisorial valuations on $R$.

Given $v\in {\rm D}(R)$ and an ideal $I\subset R$, one denotes
$v(I):=\min\{v(f)\,|\, f\in I\}$ -- this is well-defined since $I$ is finitely generated and a divisorial valuation ring
is a rank one discrete valuation ring (see \cite[Definition 6.8.9 and Theorem 9.3.2]{SHbook}).

Assume now that $R=\oplus_{_{\kern-2pt m\geq 0}} R_m$ is a standard graded domain over a field $k$.

\begin{Definition}\rm
Given a homogeneous ideal $I\subset R$ and an integer $m\geq 0$, set
$$\tilde{I}(m)=\bigcap_{v \in {\rm D}(R)}\{f\in R_m\,|\, v(f)\geq v(I)\}$$
\vskip-15pt
\noindent and $\tilde{I}=\oplus_{_{\kern-2pt m\geq 0}}\tilde{I}(m)$.
\end{Definition}

Since $v$ is a valuation, it is clear that $\tilde{I}(m)$ is a vector subspace of $R_m$ and $\tilde{I}$ is a homogeneous
ideal such that $\tilde{I}_m=\tilde{I}(m)$ for every $m\geq 0$. We will call $\tilde{I}$ the {\em divisorial cover ideal} of $I$.
We next give a more ideal theoretic formulation of this ideal:

\begin{Lemma}\label{compact_divisorial}
Notation being as above, let in addition $R_v$ denote the valuation ring on $R$ corresponding
to a given $v\in {\rm D}(R)$. Then
$$\tilde{I}=\bigoplus_{m\geq 0}\,\left(\bigcap_{v\in {\rm D}(R)}(IR_v\cap R_m)\right).$$
\end{Lemma}
\demo Note that both sides of the sought equality are homogeneous ideals of $R$, hence it suffices to show the equality
of the homogeneous parts of the same degree.
Thus, fix  $m\geq 0$ and $v\in {\rm D}(R)$.
Let $g\in R$ be such that $v(g)=v(I)$.
Then $f\in R_m$ belongs to $\tilde{I}_m$ if and only if $v(f/g)\geq 0$, i.e, if and only if $f/g\in R_v$.
This shows that $\tilde{I}_m=gR_v\cap R_m\subset IR_v\cap R_m$.
Conversely, if $f=gu\in IR_v\cap R_m$, with $g\in I$ and $u\in R_v$, then $v(f)=v(g)+v(u)\geq v(g)\geq v(I)$,
hence $f\in \tilde{I}_m$.
Since $f$ is arbitrarily taken in $IR_v\cap R_m$, it follows that $IR_v\cap R_m\subset \tilde{I}_m$.
\qed

\smallskip

Given an ideal $I\subset R$, let $\bar{I}\subset R$ denote its integral closure.

\begin{Corollary}\label{inclusion_closure_cover}
Let $R$ be a standard graded domain over a field and let $I\subset R$
denote a homogeneous ideal. Then $\bar{I}\subset \tilde{I}$.
\end{Corollary}
\demo First note that the integral closure $\bar{I}$ of $I$ is still homogenous.
On the other hand, one has $\bar I=\bigcap_{v} (IR_v\cap R)$, where this time around $v$ varies over all rank one discrete valuations of $Q$
positive on $R$ \cite[Proposition  6.8.2]{SHbook}.
Clearly, $\bar I=\bigcap_{v} (IR_v\cap R)=(\bigcap_{v} IR_v)\cap R$, hence $\bar{I}_m=(\bigcap_{v} IR_v)\cap R_m=
\bigcap_{v} (IR_v\cap R_m)$ for $v$ varying over the rank one discrete valuations  positive on $R$.
Since the divisorial valuations on $R$ are among the latter, we can apply Lemma~\ref{compact_divisorial} to conclude.
\qed

\medskip

For the next result we need the following basic result, part of which holds more generally in any dimension.
For the sake of both simplicity and objectiveness we stick to the $2$-dimensional case.

\begin{Lemma}\label{order_is_valuation}
Let $X$ stand for a smooth surface over $k$ and let $p\in X$. Denote by $(\mathcal{O}_p,\,\fm_p)$ the local ring of $p$ on $\pp^2$
and its unique maximal ideal.
Let $\mathfrak{o}_p$ stand for the order function relative to $\fm_p$.
Then:
\begin{itemize}
\item[{\rm (a)}] $\mathfrak{o}_p$ extends to  divisorial valuation $v_p$ centered on $\mathcal{O}_p$.
\item[{\rm (b)}] Let  $F:\pp^2\dasharrow\pp^2$ be a rational map with base ideal $I$ and weighted cluster $\mathcal{K}=
(K, \boldsymbol\mu)\,;$ then $v_p(I)=\mu_p$ for every $p\in K$.
\end{itemize}
\end{Lemma}
\demo
(a) This is a well-known fact going back to Zariski {\em et. al} -- a good reference is \cite[Theorem 6.7.9]{SHbook}.

(b) This part follows from the fact that any proper or infinitely near point of $K$
is a point on some (smooth) surface obtained as a successive blowup of a point on $\pp^2$ -- so one can apply part (a) --
and the definition of $\mu_p$. Indeed, tracing through this notion one sees that $\mu_p$ is obtained by evaluating
the order function at that step on the corresponding extension of the ideal $I$. To see the latter, observe that,
at each infinitely near point $p$, $\mu_p$ takes value on a nonempty open subset of the parameters.
\qed

\begin{Proposition}\label{from_divisorial_to_full}
Let $I\subset R$ be the base ideal of a rational map $F:\pp^2\dasharrow\pp^2$ with weighted cluster $\mathcal{K}=
(K, \boldsymbol\mu)$.
Then $\tilde{I}\subset I_{\mathcal{K}}$.
\end{Proposition}
\demo
For every $p\in K$ let $v_p$ as in Lemma~\ref{order_is_valuation}(a) denote the divisorial valuation ring induced
by the order function defined on the local ring $(\mathcal{O}_p,\,\fm_p)$.
Then, by Lemma~\ref{order_is_valuation},(a) and(b), one has
$$\tilde{I}_m= \kern-8pt\bigcap_{v \in {\rm D}(R)}\{f\in R_m\,|\, v(f)\geq v(I)\} \kern-3pt\subset
\bigcap_{p \in \mathcal{K}}\{f\in R_m\,|\, v_p(f)\geq v_p(I)\}=\kern-3pt\bigcap_{p \in \mathcal{K}}\{f\in R_m\,|\, v_p(f)\geq \mu_p\},$$
for every $m\geq 0$.

On the other hand, $(I_{\mathcal{K}})_m=\{f\in R_m\,|\, \mathbf{P}_{\mathcal{K}}^{-1}\cdot
 (\mathbf{e}_K(f)-{\boldsymbol\mu}_K)^t\geq 0\}$.
Now, for $p\in K$, the effective multiplicity $e_p$ is also given by the order function, hence $e_p(f)=v_p(f)$
for every $p\in K$.
As already observed and used,  $\mathbf{P}_{\mathcal{K}}^{-1}$ has only nonnegative entries.
Therefore, $\mathbf{e}_K(f) - {\boldsymbol\mu}_K=\mathbf{v}_K(f)-{\boldsymbol\mu}_K\geq 0$ certainly implies that
$\mathbf{P}_{\mathcal{K}}^{-1}\cdot (\mathbf{e}_K(f)-{\boldsymbol\mu}_K)^t\geq 0$.

This concludes the proof of the statement.
\qed

\medskip

We collect the various results in the following

\begin{Theorem}\label{almost_all_inclusions}
Let $I\subset R$ be the base ideal of a rational map $F:\pp^2\dasharrow\pp^2$ with weighted cluster $\mathcal{K}$.
Then
$$I\subset \bar{I}\subset \tilde{I}\subset I_{\mathcal{K}}\subset \mathfrak{F}_{\mathcal{K}}.$$
\end{Theorem}

\subsection{Complements}

We state a basic fact of standard graded rings, reminiscent of
an argument in the proof of \cite[Theorem 5.19]{Har}.
It gives an Artin--Rees type of equality ``on the nose'' (see \cite[Chapter 13]{SHbook} for the known
aspects relating Artin--Rees and integral closure in more generality).

 \begin{Proposition}\label{linear_multiple}
 Let $(R,R_+)$ denote a standard graded ring over a field $k$ and its maximal irrelevant ideal.
 Write $n:=\edim(R)=\dim_kR_1$.
 Let $I\subset R_+$ stand for a homogeneous ideal of codimension $\geq 2$ generated by $n$ $k$-linearly independent
 forms of degree $d$, for some $d\geq 1$.
 Given a homogeneous ideal $J\subset R_+$ such that  $I_d=J_d$, one has:
\begin{enumerate}
\item[{\rm (a)}] For any integer $j\geq 0$,
\begin{equation}\label{linear_multiple_subset}
R_+^jJ\cap I=R_+^jI.
\end{equation}
\item[{\rm (b)}] If, moreover, $J\subset I:R_+^{\infty}$ then $J$ is integral over $I$.
\end{enumerate}
 \end{Proposition}
\demo
Clearly, $I\subset J$.
Moreover, $J_{\ell}=\{0\}$ for $\ell<d$.
Indeed, if $0\neq g\in J_{d-1}$ then $R_1g\subset J_d=I_d$.
Since $R_1g$ is spanned by $k$-linearly independent elements and $I_d$ has $k$-vector dimension $n$,
it follows that $I=(R_1g)$, hence $I$ would have codimension $\leq 1$.

To prove (a), it clearly suffices to show the inclusion $R_+^jJ\cap I\subset R_+^jI$.
Now, since $R_+$ is the maximal irrelevant ideal of $R$, then $R_j=(R_+^j)_j$ for any $j\geq 0$.
Therefore, since $I_{\ell}=\{0\}$ for $\ell<d$, one can write $I_{d+j}= R_jI_d=(R_+^j)_jI_d$,
for any $j\geq 0$ (any linear form that multiplies the generators of $I$ into $I_{d+1}$ can be absorbed by
$R_+$).

Since $R_+^jJ\cap I\subset I$, it follows that
$$(R_+^jJ\cap I)_{d+j}\subset I_{d+j}=(R_+^j)_jI_d\subset (R_+^jI)_{d+j},
$$
for any $j\geq 0$.
On the other hand, we have  $J_{\ell}=\{0\}$ for $\ell<d$ as shown above.
Therefore $(R_+^jJ)_{\ell+j}=\{0\}$ for $\ell<d$ and any $j\geq 0$.
Thus, finally one has $R_+^jJ\cap I\subset R_+^jI$ for any $j\geq 0$, as required.

To prove (b), one notes that the hypothesis implies the inclusion $R_+^sJ\subset I$ for some $s\geq 0$.
Then clearly $R_+^sJ=R_+^sJ\cap I\subset R_+^sI$.
Therefore the result follows from the so-called {\em determinantal trick} (see \cite[Corollary 1.1.8]{SHbook}).
\qed

\medskip

As an immediate consequence of this criterion and the results of the previous subsection,
we file the following result.

 \begin{Proposition}\label{plane_pan_russo}
Let $I$ be the base ideal of a Cremona map of degree $d\geq 1$ of $\pp^2$ with weighted cluster $\mathcal{K}$.
Then
$$I\subset I^{\rm sat}\subset\bar{I}\subset  \tilde{I}\subset I_{\mathcal{K}}\subset \mathfrak{F}_{\mathcal{K}}.$$
{\sc Supplement.}  $I_d = I^{\rm sat}_d = \bar{I}_d =  \tilde{I}_d = (I_{\mathcal{K}})_d$.
 \end{Proposition}
 \demo
By  \cite[Proposition 1.2]{PanRusso}, $I_d=I^{\rm sat}_d$.
Thus, we can apply Proposition~\ref{linear_multiple}(b) to conclude that $I^{\rm sat}\subset\bar{I}$.
The other inclusions were proved in Theorem~\ref{almost_all_inclusions}.

The statement in the supplement follows from \cite[Proposition 2.5.2]{alberich} in which it is shown that
$I_d=(I_{\mathcal{K}})_d$.
\qed

\medskip

Note that \cite[Proposition 1.2]{PanRusso} translates into the
basic assumption of the previous section, namely, that $\beg (I^{\rm sat}/I)\geq d+1$.
Thus, Theorem~\ref{ldual} and all of its consequences are immediately applicable and will
subsequently be drawn upon.
For example, one has:

\begin{Proposition}\label{cupperbound}
Let $I=(I_d)\subset R=k[x_0,x_1,x_2]\,(d\geq 2)$ stand for the base ideal of a plane Cremona map without fixed part.
If $I$ is non-saturated then $d\geq 5$ and:
\begin{itemize}
\item[{\rm (i)}]
$(d{\,}^2+3d-4)/2\leq e(R/I)\leq (5d\,^2-21d)/2$, where lower bound is for the case where the base points are proper,
while the upper bound is valid for $d\geq 6$
\item[{\rm (ii)}]
$2d-4-\lfloor(d-5)/2\rfloor\leq \Reg(R/I)\leq 2d-4$.
\end{itemize}
\end{Proposition}
\demo
The bounds for the regularity follows from Theorem~\ref{from_lower_bound} (ii) and  the inequalities
$2(\Reg(R/I)+2)\geq d_r+d_{r-1}\geq 1+3d$ in (\ref{ineq}). As to the upper bound for the multiplicity, it comes out of the following
 equality derived from (\ref{hseries})
\begin{equation}\label{multiplicity}
e(R/I)=\sum_{i=1}^r {{d_i}\choose {2}}- \sum_{m=1}^{r-2} {{D_m}\choose {2}} -3{{d}\choose {2}}.
\end{equation}
The lower bound $e(R/I)\geq e(R/\mathfrak{F}_{\mathcal{K}})=(d^2+3d-4)/2$
has been previously explained, where $\mathfrak{F}_{\mathcal{K}}$ is the fat ideal as above
and after using the equations of condition (\ref{eqs_condition}).
\qed

\begin{Remark}\rm It would be interesting to obtain an estimate of $e(R/I_{\mathcal{K}})$ in order
to get tighter lower bounds for $e(R/I)$.
\end{Remark}

Observe that, in contrast to the intangible ideal theoretic properties of the full ideal of curves $I_{\mathcal{K}}$,
the fat ideal based on the proper part of $\mathcal{K}$ is unmixed and integrally closed.
An immediate consequence of this facet is the following.

\begin{Corollary}\label{sat_is_closure}
Let $F$ be a plane rational map without fixed part and base ideal $I$, and let $\mathfrak{F}_{\mathcal{K}}$
 denote the fat ideal associated to the
proper part of the corresponding weighted cluster of $F$.
If $\mathfrak{F}_{\mathcal{K}}=I^{\rm sat}$ then $\mathfrak{F}_{\mathcal{K}}$ is the integral closure of $I$.
\end{Corollary}
\demo
Since $\mathfrak{F}_{\mathcal{K}}$ is integrally closed and $I\subset \mathfrak{F}_{\mathcal{K}}$ then the integral closure
$\bar{I}$ of $I$ is contained in $\mathfrak{F}_{\mathcal{K}}$.
By Proposition~\ref{linear_multiple} (b), $\mathfrak{F}_{\mathcal{K}}\subset \bar{I}$. Therefore, $\mathfrak{F}_{\mathcal{K}}=\bar{I}$.
\qed
\begin{Example}\rm
An instance of the above is a simple Cremona map of homaloidal type $(5\,; 2^6)$ (see the proof of
Theorem~\ref{proper_types_degree5}).
\end{Example}

Finally, using
Proposition~\ref{plane_pan_russo} and Theorem~\ref{from_lower_bound} one gets:

\begin{Corollary}\label{degree_at_most_4}
Let $I=(I_d)\subset R=k[x_0,x_1,x_2]\,(d\geq 2)$ stand for the base ideal of a plane Cremona map without fixed part.
If $d\leq 4$ then $I$ is saturated.
\end{Corollary}

\section{Steps in the classification of plane Cremona maps}

\subsection{de Jonqui\`eres maps}

Among plane Cremona maps, the so-called {\em de Jonqui\`eres map} plays a fundamental role
going back at least to Castelnuovo's celebrated proof of Noether theorem (\cite[Proposition 8.3.4]{alberich}.
Following \cite[2.6.10]{alberich}, we define it  as a plane Cremona map of degree $d\geq 2$ whose homaloidal type is
$(d\,; d-1,1^{2d-2})$.

A basic geometric datum is its close association with the so-called {\em monoids}.
Besides, this class of maps enjoys many interesting properties from the algebraic and homological
viewpoints.
It may be convenient to  give an overview of some of these properties and describe some
relevant families of such maps.

\subsubsection*{Monomial de Jonqui\`eres maps}

An important class of Cremona maps is that of monomial Cremona maps.

\begin{Definition}\rm
A {\em monomial Cremona map} is a Cremona map in $\pp^n$ whose base ideal is generated by monomials
in $k[x_0,\ldots, x_n]$.
\end{Definition}

There is as of now a reasonably extensive literature on these maps (see, e.g.,
\cite{Barbara}, \cite{CoSi}, \cite{Andre}, \cite{SiVi}, \cite{SimisVilla}, \cite{CremonaMexico}).
It is on itself a guiding case study. The following proposition covers some
basic properties of plane such maps and gives a characterization
of plane monomial de Jonqui\`eres  maps; parts (c) and (d) below seem to be new.

\begin{Proposition}\label{monomial_jonquieres}
Let $F:\pp^2\dasharrow \pp^2$ denote a plane monomial Cremona map
and let $I\subset R$ stand for its base ideal.
Then:
\begin{enumerate}
\item[{\rm (a)}] Up to permutation of the variables {\rm (}source{\rm )} and the defining
 monomials {\rm (}target{\rm )}, the base ideal $I$ is one of the following:
\begin{itemize}
\item $(xy,xz,yz)$
\item $(x^d,x^{d-1}y,y^{d-1}z)$, with $d\geq 1$
\item $(x^d, x^{d-(a+b)}y^az^b, y^{d-c}z^c)$, with $abc\neq 0$, $d\geq a+b$ and $ac-b(d-c)=\pm 1$.
\end{itemize}
\item[{\rm (b)}] $I$ is saturated.
\item[{\rm (c)}] $F$ is a de Jonqui\`eres map exactly in the following cases:
\begin{itemize}
\item $(xy,xz,yz)$
\item $(x^d,x^{d-1}y,y^{d-1}z)$, with $d\geq 2$
\item $(x^d, xyz^{d-2}, yz^{d-1})$, with $d\geq 3$.
\end{itemize}
\item[{\rm (d)}] If $F$ is a de Jonqui\`eres map then $I$ is an integrally closed ideal exactly in the following cases:
\begin{itemize}
\item $(xy,xz,yz)$
\item $(x^2,xy,yz)$
\item $(x^3,x^2y,y^2z)$
\item $(x^3,xyz,yz^2)$.
\end{itemize}
\end{enumerate}
\end{Proposition}
\demo
(a) This part is proved in \cite[Lemma 3.3 and Proposition 3.5]{CremonaMexico}.

\smallskip

(b) This is a result of \cite[Proposition 2.11]{Andre} drawing upon part (a).

\smallskip

(c)  We first argue that the three alternatives give in fact
de Jonqui\`eres maps according to the above definition.
We skip the discussion of the first case as being sufficiently known.

It suffices to show
that there is a proper base point  where the defining $d$-forms
have minimum multiplicity (i.e., order) $d-1$, for then the equations of condition
yield the existence of $2d-2$ additional  proper or infinitely
near base points.

Consider the case $I=(x^d,x^{d-1}y,y^{d-1}z)$, $d\geq 2$.
The map has only two proper base points: $p=(0:0:1)$ and $q=(0:1:0)$.
By passing to the respective local rings, one readily finds that the minimum
multiplicity at $p$ is $d-1$, while at $q$ it is $1$.

In the case $I=(x^d, xyz^{d-2}, yz^{d-1})$, we have again the same two proper base points.
Again, an immediate check gives multiplicity $d-1$ at the point $q$ (and $1$ at $p$).

\smallskip

Conversely, let $F$ be a de Jonqui\`eres map.
By (a), we may assume that its base ideal has the third form.
Because of the restrictions on the integers $a,b,c$ in this case, we again find that the map
has only the two proper base points $p,q$ as before.
Therefore, one (and only) one of these points has minimum multiplicity $d-1$
for the curves of the system.

Now, locally at $q$ and $p$ the minimum multiplicity of the curves of the system is
$$\min\{d-a,c\}\quad  {\rm and}\quad \min\{d-(a+b)+a=d-b,d-c\},$$
respectively.
But since permuting $y$ and $z$ will not change the form of the map,
we may assume that this point is $q=(0:1:0)$.
Suppose first that $d-1=d-a$, hence $a=1$.
Then, on the other point the minimum multiplicity must be $1$, hence either $b=d-1$
or $c=d-1$. If $b=d-1$ then $a+b=d$, contradicting one of the restrictions in
the third case of part (a).
Therefore, it must be the case that $c=d-1$. But then $ac-b(d-c)=\pm 1$
now reads $b=d-1\pm 1$. Since $b<d$, we must have $b=d-2$, as required.

The discussion of the alternative $c=d-1$ at the outset is similar: $ac-b(d-c)=\pm 1$
now gives $a\leq (a+b)a-b\leq (d-1)a-b=\pm 1$. Since $a>0$, we get  $a=1$.
Then as before, $b=d-2$.

\smallskip

(d)
It is well-known that $(xy,xz,yz)$ is even a normal ideal (see, e.g., \cite[Theorem 1.1 and Corollary 2.8]{bowtie},
also \cite{HiOh}) and so is $(x^2,xy,yz)$ by a similar token: the defining ideal of the Rees algebra in either case
is a codimension $2$ complete intersection whose Jacobian ideal has codimension $2$.

The cases of $(x^3,x^2y,y^2z)$ and $(x^3,xyz,yz^2)$ can be dealt with by readily verifying that the ideal on the local pieces
$z=1$ and $y=1$ is integrally closed; from this, an easy checking shows that the original ideal is integrally closed.

On the other hand, neither $I=(x^d,x^{d-1}y, y^{d-1}z)$ nor $J=(x^d, xyz^{d-2}, yz^{d-1})$ is integrally closed for $d\geq 4$.
Indeed,  for example:
$$\left\{
\begin{array}{lll}
x^{d/2}y^{d/2}z\notin I, & \mbox{but $(x^{d/2}y^{d/2}z)^2\in I^2,$} & \mbox{if $d\geq 4$ is even}\\[5pt]
x^{(d-1)/2}y^{(d+1)/2}z\notin I, & \mbox{but $(x^{(d-1)/2}y^{(d+1)/2}z)^2\in I^2,$} & \mbox{if $d\geq 5$ is odd}
\end{array}
\right.
$$
and, similarly
$$\left\{
\begin{array}{lll}
x^{d-1}y^2z^{(d-2)/2}\notin J, & \mbox{but $(x^{d-1}y^2z^{(d-2)/2})^2\in J^2,$} & \mbox{if $d\geq 4$ is even}\\[5pt]
x^{d-1}y^{2}z^{(d-1)/2}\notin J, & \mbox{but $(x^{d-1}y^{2}z^{(d-1)/2})^2\in J^2,$} & \mbox{if $d\geq 5$ is odd}.
\end{array}
\right.
$$
\qed

\medskip

\subsubsection*{Arbitrary de Jonqui\`eres maps}

The next proposition states the equivalence between two notions of a de Jonqui\`eres map.
Although this equivalence is used liberally (see, e.g.,  \cite[Lemme 2]{PanBoletim}, \cite[Exemple 1.3]{Pan})
we could not trace through the literature a fully rigorous algebraic proof.
Since the definition we employ here is based on properties of the weighted cluster of the map,
while the second notion essentially deals with the format of the base ideal, it would seem desirable to have such
a precise proof.

In addition, there is a very good reason to be able to navigate between the two notions as, at one end, the first is computable in terms of
any set of generators of the base ideal, while the second one is impossible to verify as its formulation
depends upon a projective coordinate change both in the source and the target of the map.
At the other end, the second gives a good handling of the algebraic properties of the ideal, whereas the first notion
falls behind in this regard.

A $z$-{\em monoid} is a $d$-form $f_{d-1}z+f_d\in R$, where $f_{d-1},f_d$ are forms in $k[x,y]$
of respective degrees $d-1,d$.
It is noteworthy  that such a form is irreducible if (and only if) $\gcd\{f_{d-1},f_d\}=1$.

\begin{Proposition}\label{characterization_of_jonquieres}
Let $F:\pp^2\dasharrow \pp^2$ be a rational map of degree $d$ with no fixed part and let
$I\subset R$ denote its base ideal.
The following conditions are equivalent:
\begin{enumerate}
\item[{\rm (i)}] $F$ is a de Jonqui\`eres map
\item[{\rm (ii)}] Up to permutation of the variables {\rm (}source{\rm )} and the defining
forms {\rm (}target{\rm )}, the base ideal $I$ is generated by $d$-forms $\{f, xq, yq\}$ such that
$f$ and $q$ are both $z$-monoids and $f$ is irreducible.
\end{enumerate}
\end{Proposition}
\demo
(i) $\Rightarrow$ (ii)
Since there must be at least one proper base point, we may assume that $p=(0:0:1)\in \pp^2$ is the
base point with the (unique) prescribed multiplicity $d-1$.
On the other hand, by the genus formula,  any irreducible $d$-form having multiplicity $d-1$ at $p$
has no other  (proper or infinitely near) singularities.
But among these forms, there is a $z$-monoid, namely, $f=zf_{d-1}(x,y)+f_d(x,y)$, with $f_{d-1}$ and $f_d$
forms of degrees $d-1$ and $d$, respectively, with $\gcd(f_{d-1},f_d)=1$.
Now, since $f_{d-1}$ and $f_d$ can be chosen generically, we may assume that they have been chosen so that,
moreover, $f$ goes through the remaining (proper or infinitely near) base points of $F$.
It follows from \cite[Proposition 2.5.2]{alberich} that $f\in I_d$.
Such a form being irreducible,  the corresponding homogeneous coordinate ring $k[x,y,z]/(f)$
is an integral domain.
Also, $V(f)$ is the image of the rational map
$$(s_0f_{d-1}(s_0,s_1): s_1f_{d-1}(s_0,s_1):-f_d(s_0,s_1)) :\pp^1_{s_0,s_1}\dasharrow \pp^2_{x,y,z},$$
which is birational onto $V(f)$.
Therefore, the natural inclusion
$$k[s_0f_{d-1}(s_0,s_1), s_1f_{d-1}(s_0,s_1), -f_d(s_0,s_1]\subset k[s_0,s_1]^{(d)}$$
into the $d$th Veronese is an equality at the level of the respective fraction fields, which in turn induces an injective $k$-homomorphism
$$k[x,y,z]/(f)\stackrel{\alpha}{\simeq} k[s_0f_{d-1}(s_0,s_1), s_1f_{d-1}(s_0,s_1), -f_d(s_0,s_1]\subset k(s_0,s_1).
$$

Extending $f$  to a $k$-vector base $\{g,h,f\}$ of $I_d$ corresponds to a coordinate change on the
target, hence we may assume that the map
$F:\pp^2\dasharrow \pp^2$ is defined by $(g:h:f)$.
On the other hand, for general choices of $f_{d-1}$ and $f_d$,
 the restriction of $F$ to $V(f)\subset \pp^2$ is birational onto the image.
 This gives, again and over, that the injection
 $$\sigma: k[g,h,f]/k[g,h,f]\cap (f)\injects (k[x,y,z]/(f))^{(d)}$$
 extends to an equality of the respective fields of fractions.
 But the natural map
 $$k[g,h]\stackrel{\gamma}{\injects} k[g,h,f]/k[g,h,f]\cap (f)$$
is injective (which is the translation of the fact that the restriction of $F$ maps $V(f)$ onto a
 line).
Altogether, one has an injective $k$-homomorphism
$\beta: k[g,h]\injects k(s_0,s_1)$, $\beta=\alpha\sigma\gamma$,
which is an isomorphism at the level of $d$th Veronese fields of fractions.
Thus, $(\beta(g):\beta(h))$ defines a Cremona map of $\pp^1$, hence must be equivalent to the identity map of $\pp^1$.
In other words, $\beta(g)$ and $\beta(h)$ are multiples of $s_0,s_1$, respectively, with same common factor in $k[s_0,s_1]$.
Extending $\alpha$ to the level of fractions modulo $(f)$ yields
$$\alpha(y/x)=s_1f_{d-1}(s_0,s_1)/s_0f_{d-1}(s_0,s_1)=s_1/s_0=\alpha(h/g),$$
thus implying $yg\equiv xh \pmod {f} $. Therefore,
we obtain a relation $yg-xh+\ell f=0$, for some linear form $\ell\in R$.
On the other hand, we can write $f=xp_1+yp_2$ for suitable $(d-1)$-forms $p_1,p_2\in k[x,y,z]$.
Substituting in the linear relation for $\ell f$, and using that $\{x,y\}$ is a regular sequence,
yields $g=xq-\ell p_2, h=yq+\ell p_1$, for some $(d-1)$-form $q$.
Then, as one readily checks, $f,g,h$ are (up to a sign) the $2$-minors of the matrix
$$\phi=
\left(
\begin{array}{cc}
x& p_2\\
-y& p_1\\
\ell & q
\end{array}
\right).
$$
Assume first that $d\geq 3$.
We claim that $x,y,\ell$ are linearly dependent.

Indeed, otherwise the ideal $I_1(\phi)$ would have
codimension $3$, thus implying that $I$ is generically a complete intersection.
But, since $I$ is an almost complete intersection, it would be an ideal of linear type (see, e.g.,
\cite[Proposition 3.7]{SV}).
However, since $F$ is a Cremona map and , this would contradict the main birationality criterion of \cite{AHA}.

This means that $\ell$ is a $1$-form in $k[x,y]$, hence up to a change of coordinates we may assume that
$\ell=0$. It follows that $I$ is now generated by $\{f, xq, yq\}$.
Here $f$ is an irreducible $z$-monoid by construction, while $q$ is again a $z$-monoid because
$xq$ and $yq$ have multiplicity $\geq d-1$ at $(0:0:1)$.

In the case where $d=2$, up to a coordinate change on the source and the target
 the respective base ideal has one of the following set of generators:
\begin{enumerate}
\item $F$ is given by $\{xy, xz, yz\}$ (Three  distinct proper base points $p_1,p_2,p_3$)
\item $F$ is given by $\{x^2, xy, yz\}$ (Two distinct proper base points $p_1,p_2$ and one infinitely near
point to $p_1$)
\item $F$ is given by $\{x^2, xy, y^2-xz\}$ (A unique proper base point $p$ and two infinitely near
points to $p$).
\end{enumerate}
We deal with the first of these alternatives as the other two are treated similarly.
Namely, we take the following modified set of generators: $f=xz+yz+xy=(x+y)z+xy= x(y+z)+yz, g=xz,h=yz$.
Then $f$ is an irreducible $z$-monoid, while $q=z$.
The corresponding syzygy matrix is
$$\phi=
\left(
\begin{array}{cc}
x&z\\[3pt]
-y& y+z\\[3pt]
0 & z
\end{array}
\right).
$$

\smallskip

(ii) $\Rightarrow$ (i)
According to \cite[Lemme 2]{PanBoletim}, the rational map with these properties is birational.
In particular, the equations of condition are satisfied. Thus it suffices to show that
the minimum multiplicity at the base point $(0:0:1)$ is $d-1$.
But this follows immediately from the monoidal structure of $f$ and $q$.
\qed

\begin{Remark}\rm
A self-contained argument for proving that the map as given in (ii) is birational comes from the details
in the proof of Theorem~\ref{Rees_of_Jonq}(i). Indeed, in [loc. cit.] it is shown in particular that there
exists a relation of $\{f,xq,yq\}$ which is a form of bidegree $(1,d-2)$.
Such a form, along with the form of bidegree $(1,1)$ coming from the trivial relation of $\{xq,yq\}$,
implies that the pertinent matrix in \cite[Theorem 1.18]{AHA} has rank $2$, thus forcing birationality.
\end{Remark}

Recall that an ideal $I$ of a ring $R$ is {\em of linear type} if the natural surjective
homomorphism from the symmetric algebra of $I$ to its Rees algebra is injective.

The following ideal theoretic properties follow suit.

\begin{Corollary}\label{jonquieres}
Let $F$ denote a plane de Jonqui\`eres map of degree $d\geq 2$ with base ideal $I\subset R=k[x,y,z]$.
Then:
\begin{itemize}
\item[{\rm (a)}]  $R/I$ is Cohen--Macaulay, hence $I=I^{\rm sat}$.
\item[{\rm (b)}] The degree {\rm (}algebraic multiplicity{\rm )} of the scheme $R/I$ is $d(d-1)+1$.
\item[{\rm (c)}] $I$ is an ideal of  linear type if and only if $d=2$.
\end{itemize}
\end{Corollary}
\demo
(a)
Then a set of generators of $I$  are the $2\times 2$ minors of the matrix
$$\phi=
\left(
\begin{array}{cc}
x&p_2\\
-y& p_1\\
0 & q
\end{array}
\right).
$$
Since $I$ has codimension $\geq 2$, this shows that $R/I$ is a Cohen--Macaulay ring, hence $I$ is an unmixed ideal.
Thus, a de Jonqui\`eres map has a saturated base ideal.

\smallskip

(b) It follows from (a) that $R/I$ has a minimal graded resolution of the form
\begin{equation}\label{res_of_jonquieres}
0\rar R(-(2d-1))\oplus R(-(d+1))\lar  R(-d)^3\lar R.
\end{equation}
Therefore, its Hilbert series is $(1-3t^d+t^{d+1}+t^{2d-1})/(1-t)^3$.
From this, by taking the second $t$-derivative of the numerator and evaluating at $t=1$,
one obtains the degree
$$e(R/I):=\frac{-3d(d-1)+(d+1)d+ (2d-1)(2d-2)}{2}=d(d-1)+1.$$

\smallskip

(c) As mentioned in the proof of Proposition~\ref{characterization_of_jonquieres}, since $I$ is an almost
complete intersection one knows that $I$ is of linear type if and only
it is generically a complete intersection,
i.e., if and only if the ideal of $1$-minors of its
structural Hilbert--Burch matrix $\phi$ as above is $(x,y,z)$-primary. By the
monoidal form of $q$,
$q$  has a pure power term in $z$ only if $\deg(q)=1$; else, it is either
$p_1$ or $p_2$ that has a nonzero term in $z^r$ alone, for some $r\geq 1$.
Again, from the monoidal form of $f$ above this happens if and only if $\deg (p_1)=1$.
\qed

\begin{Remark}\rm
(1) The interest of a calculation such as the one in item (b) above is that, even if the points are proper
and in general position, the corresponding $\boldsymbol\mu$-fat ideal $\mathfrak{F}_{\mathcal{K}}$
 may not coincide with $I^{\rm sat}$, where $I$ denotes the base ideal of the linear system
in degree $d$. An example is obtained with $d=4$, where $I$ is saturated by
Corollary~\ref{jonquieres}(a). One has $e(R/I)=13$ by the above corollary, while
$e(R/\mathfrak{F}_{\mathcal{K}})=12$ by (\ref{multplicity_of_fat}).

(2) As embodied in the proof of Proposition~\ref{characterization_of_jonquieres}, the statement of (c) in the
above proposition is also a consequence of the birationality criterion
of \cite{AHA}.
\end{Remark}

\medskip

One can even go one step further to give the structure of the Rees algebra of the base ideal of
a de Jonqui\`eres map, in fact, of any rational map $\pp^2\dasharrow \pp^2$ whose base ideal has a similar structure.
More precisely, we prove:

\begin{Theorem}\label{Rees_of_Jonq}
Let $I\subset R$ denote the base ideal of a rational map $\pp^2\dasharrow \pp^2$  of degree $d\geq 2$
with no fixed part, whose  syzygy matrix has the form
$$\phi=
\left(
\begin{array}{cc}
x&p_2\\
-y& p_1\\
0 & q
\end{array}
\right),
$$
where $q$ has a nonzero term which is a pure power of $z$ if and only if $d=2$.
Let $\mathcal{R}(I)$ denote the Rees algebra of $I$.
One has:
\begin{enumerate}
\item[{\rm (i)}] If $\mathcal{R}(I)\simeq R[t,u,v]/\mathcal{J}$ stands for a minimal presentation
then $\mathcal{J}$ is minimally generated by $d$
polynomials of bidegrees $(1,1), (d-1,1), (d-2,2),...,(1,d-1)$
\item[{\rm (ii)}] $I^j$ is saturated for $1\leq j\leq d-1$, while $I^d$ is not saturated
\item[{\rm (iii)}] $\mathcal{R}(I)$ is Cohen--Macaulay if and only if $d\leq 3$.
\end{enumerate}
\end{Theorem}
\demo
(i) For $d=2$ we are assuming that $q=z+q'$, where $q'\in k[x,y]$.
Therefore, the ideal $I_1(\phi)=(x,y,z)$. This forces $I$ to be generically a complete intersection,
hence is of linear type -- an argument we have repeatedly used in the text.
Thus, the result is immediate in this degree.

\smallskip

Suppose that $d\geq 3$. The hypothesis on the form of $\phi$ implies at the outset that $I_1(\phi)=(x,y)$.
This is suited to applying the method of the {\em Sylvester forms} as indicated in \cite{syl1}
and based on the standard bigrading of the polynomial ring $S=R[t,u,v]$, where $x,y,z$ have bidegree
$(1,0)$ and $t,u,v$ have bidegree $(0,1)$.

 The Sylvester form associated to a set of polynomials
 is computed with respect to a given ideal that serves as a kind of ``frame'' for these polynomials.
 The advantage of this procedure is that the forms obtained are among the defining Rees equations of $I$.

Starting with the Rees equations $F=xt-yu, G=p_2t+p_1u+qv\in \mathcal{J}$,
 coming from the syzygies of $I$, we can then write

 \[\left[\begin{array}{c} F \\ G \end{array}\right]=
                     \left(
    \begin{array}{cc}
      t & -u \\
      G_x & G_y\\
    \end{array}
  \right)
                      \left[\begin{array}{c} x
                      \\ y \end{array} \right],
\]
for suitable forms $G_x,G_y\in R[t,u,v]$ of bidegree $(d-2,1)$.
By Cramer's rule, $H_1:=\det(C)\in \mathcal{J}$, where
$$C=
\left(
    \begin{array}{cc}
      t & -u \\
      G_x & G_y\\
    \end{array}
  \right).
$$
This is our first Sylvester form.
We have $H_1\neq 0\,$; indeed, otherwise $[P \;Q ]\cdot C=0$, for some nonzero row vector with $P,Q\in R[t,u,v]$.
But then the transpose $[P \; Q]^t$ would be a syzygy of $\{F,G\}$. Since the latter is a regular
sequence, it would give that the columns of $C$ are syzygies of $\{-G,F\}$, which is absurd.

Note that $H_1$ is of bidegree $(d-2,2)$.
By the same procedure, we obtain the next Sylvester form:

 \[\left[\begin{array}{c} F \\ H_1 \end{array}\right]=
                     \left(
    \begin{array}{cc}
      t & -u  \\
      (H_1)_x & (H_1)_y\\
    \end{array}
  \right)
                      \left[\begin{array}{c} x
                      \\ y \end{array} \right],
  \]
for suitable forms $(H_1)_x, (H_1)_y\in R[t,u,v]$ of bidegree $(d-3,2)$.
By the same argument, the determinant $H_2$ of the above  $2\times 2$ ``content'' matrix is a nonzero Sylvester form
of bidegree $(d-3, 3)$ belonging to $\mathcal{J}$.
Continuing this way, always using $F$ and the updated $H_{i-1}$, we obtain a subideal
$(F,G,H_1, \ldots, H_{d-2})\subset \mathcal{J}$,
where $H_i$ is a form of bidegree $(d-(i+1), i+1)$, for $i=1,\ldots, d-2$.

Now, $H_{d-2}\notin (x,y)R[t,u,v]$ because otherwise we could construct a nonvanishing  Sylv\-ester form of bidegree
$(0,d)$. This would give a nonzero form in $\mathcal{J}\cap k[t,u,v]$, contradicting birationality.
But then the $R$-content ideal of the forms $\{F,G,H_1, \ldots, H_{d-2}\}$
has codimension $3$. By \cite{syl3}, one knows that $(F,G,H_1, \ldots, H_{d-2})=\mathcal{J}$,
as was to be shown.

\smallskip

(ii) While in (i) we stressed the bigraded structure of the Rees algebra as induced from the
bigrading of the polynomial ring $S=R[t,u,v]$, for the purpose of this part we
preliminarily focus on the usual standard $\NN$-grading as induced from the standard $\NN$-grading of
$S$ with $S_0=R$.
More explicitly, let $W$ be an additional (tag) variable of bidegree $(0,1)$.
Write $\mathcal{R}(I)=\oplus_{i\geq 0}\, I^iW^i\subset R[W]$,
so that in the minimal graded presentation $\mathcal{R}(I)\simeq S/\mathcal{J}$ the
component $I^iW^i$ of $\NN$-degree $i$
is presented as $I^iW^i\simeq S_i/\mathcal{J}_i$.
Thus, we get an exact sequence of $\NN$-graded $R$-modules
\begin{equation}\label{power_presentation} 0\rar \mathcal{J}_i\lar S_i\lar I^iW^i\rar 0.
\end{equation}

We now prove that, for $i\leq d-1$, the $R$-module $\mathcal{J}_i$ is  free. For this purpose, drawing upon
the notation in part (i), we contend that $\mathcal{J}_i$ is generated as a graded $R$-module by the elements
$$\{(t,u,v)_{i-1}F,v^{i-1}G,v^{i-2}H_1,\cdots,vH_{i-2},H_{i-1}\},$$
where $(t,u,v)_{i-1}$ denotes the set monomials of degree $i-1$ in $k[t,u,v]$.
Since the cardinality of this set is ${{i+1}\choose {2}}+i={{i+2}\choose {2}} -1$ which is
the rank of $J_i$, we will be done.

As a slight check, the elements listed above all have bidegree $(j, i)$ for various $j$'s, hence
they certainly belong to $\mathcal{J}_i$.
We proceed by induction on $i$.

For $i=1$ the list consists of $\{F, G\}$; since $\mathcal{J}_1$ always generates the presentation
ideal of the symmetric of $I$, the result follows immediately.

Consider the $2\times d$ matrix
$$\mathfrak{H}=
\left(
 \begin{array}{cccccc}
 u & -x & G_y & (H_1)_y & \hdots & (H_{d-3})_y\\
 -t & y & G_x & (H_1)_x & \hdots & (H_{d-3})_x
 \end{array}
 \right)
 $$
Then, by construction, the forms $\{F,G,H_1, \ldots, H_{d-2}\}$ are the $2$-minors of the
initial $2\times 3$ submatrix and, in addition, the remaining minors fixing the first column
(taking the ideal of all $2$-minors fixing the first two columns only gives repetitions).

We do the case $i=2$ to illustrate the pattern in the general inductive step.
Thus,  according to part (i), $\mathcal{J}_2$ is generated as a graded $R$-module by the set
\begin{equation}\label{apparent_gens}
\{(t,u,v)_1\mathcal{J}_1,H_1\}=\{tF,uF,vF,tG,uG,vG,H_1\}.
\end{equation}
Now, the forms $F,H_1, G$ are the $2$-minors (up to a sign) of the following $3\times 2$ initial submatrix of $\mathfrak{H}$:
 $$\left(
 \begin{array}{ccc}
 u & -x & G_y\\
 -t & y & G_x
 \end{array}
 \right).
  $$
Using that the two rows are relations of the minors and noting that $G_x$ and $G_y$
have bidegree $(d-2,1)$, one readily sees that
$tG, uG$ both belong to the graded $R$-submodule $R_1H_1+R_{d-2} (t,u,v)_1F$.
 Therefore in the set (\ref{apparent_gens})
of generators of $\mathcal{J}_2$, $tG$ and $uG$ are redundant.

Let now $3\leq i\leq d-1$ and  assume that $\mathcal{J}_{i-1}$  is  generated as a graded $R$-module by
$$L_{i-1}:=\{(t,u,v)_{i-2}F,v^{i-2}G,v^{i-3}H_1,\cdots,vH_{i-3},H_{i-2}\}.$$
Then, again from  part (i), $\mathcal{J}_i$ is generated as a graded $R$-module by $(t,u,v)_1L_{i-1}$ and by $H_{i-1}$.
Using this time around the submatrix
$$\left(
 \begin{array}{ccc}
 u & -x & (H_{i-2})_y\\
 -t & y & (H_{i-2})_x
 \end{array}
 \right)
  $$
 of $\mathfrak{H}$, whose $2$-minors are $F,H_{i-1}, H_{i-2}$, a similar discussion as in the case $i=2$
 shows that $tH_{i-2}$ and $uH_{i-2}$ belong to the graded $R$-submodule generated by
$H_{i-1}$ and $(t,u,v)_{i-1}F$.
Thus, the inductive step yields that, for $j=1,\cdots,i-2$, the elements
 $uv^jH_{i-j-2}$ and $tv^jH_{i-j-2}$ belong to the graded $R$-submodule generated by $v^jH_{i-j-1}$ and $v^j(t,u,v)_{i-j-1}F$,
 where $H_0:=G$, and hence are
 superfluous generators. Therefore $$\{(t,u,v)_{i-1}F,v^{i-1}G,v^{i-2}H_1,\cdots,vH_{i-2},H_{i-1}\}$$  generates
 the graded $R$-module $\mathcal{J}_i$, as needed to be shown.

To get the explicit minimal free graded $R$-resolution of $I^i$, with $i\leq d-1$, one uses the identification
$S_i=R[t,u,v]_i\simeq R\otimes_k k[t,v,u]_i\simeq R\otimes_k k^{{i+2}\choose {2}}\simeq R^{{i+2}\choose {2}}$
of graded $R$-modules and the appropriate shift based on the bidegree $(di,0)+(0,i)=(di,i)$ of
$I^iW^i$, in order to obtain
\begin{eqnarray}\label{res_of_power} \nonumber
0 &\lar & \left(R(-(d i+1))^{{i+1}\choose {2}}\right)\bigoplus \left(\oplus _{l=1}^{i} R(-(d-l+d i))\right)
\stackrel{\Psi_i}{\lar} R(-di)^{{i+2}\choose {2}}\lar I^i\rar 0,
\end{eqnarray}
where
$$\Psi_i=
\left(
\begin{array}{ll@{\vrule\quad}r}  
\hbox{$\mathfrak{L_i}$}&&\hbox{$\mathfrak{D_i}$}\\[3pt]
\end{array}
\right),
\quad\mathfrak{L_i}=
\left(
    \begin{array}{ccccc}
      0 & 0 & \cdots & 0 & x \\
      0 & 0 & \cdots & 0 & -y \\
       0 & 0 & \cdots & x &0 \\
      0 & 0 & \cdots & -y & 0\\
       \vdots & \vdots & \vdots & \vdots &\vdots\\
       x & 0 & \cdots & 0 & 0 \\

      -y & 0 & \cdots & 0 &0\\
      0 & 0 & \cdots & 0 & 0
    \end{array}
  \right)
$$
with:
\begin{itemize}
\item $\mathfrak{L_i}$ is of size $({{i+2}\choose {2}}-1)\times {{i+1}\choose {2}}$, induced by the linear syzygy of $I$;
\item $\mathfrak{D_i}$ is of size ${{i+2}\choose {2}}\times (i-1)$, induced by the syzygy of  $I$ of standard
degree $d-1$ and by the subsequent syzygies of $I^l$,  one for each Sylvester form $H_l$, $l=2,\ldots, i+1$.
\end{itemize}

\smallskip

Finally, we deal with the step $i=d$. By the previous argument $\mathcal{J}_{d-1}$ is minimally generated by  $$\{(t,u,v)_{d-2}F,v^{d-2}G,v^{d-3}H_1,\cdots,vH_{d-3},H_{d-2}\}.$$
We prove that $\mathcal{J}_d$ is minimally generated by
$$\{(t,u,v)_{d-1}F,v^{d-1}G,v^{d-2}H_1,\cdots,v^2H_{d-3},vH_{d-2},tH_{d-2},uH_{d-2}\}.$$
The elements of this list, except for the last three ones, are part of a minimal set of generators as it follows
again by a similar prior argument.
As to the last three elements $vH_{d-2},tH_{d-2},uH_{d-2}$ it is clear that they do not belong
 to the ideal generated by the previous ones since $H_{d-2}\not\in (x,y)S$ (as has been remarked in the proof of part (i)).

\medskip

(iii) If $d\leq 3$ we may assume that $d=3$ since $d=2$ is even a complete intersection.
Then $\mathcal{J}$ is generated by the $2$-minors of the matrix $\mathfrak{H}$ as above.
Clearly, this is Cohen--Macaulay.
Conversely, suppose that $\mathcal{R}(I)$ is Cohen--Macaulay.
Since the codimension of $\mathcal{J}$ is $2$, it is defined by a Hilbert--Burch type of matrix
again. But the presence of a minimal generator of bidegree $(1,1)$ forces this matrix to be $3\times 2$.
Therefore,  $\mathcal{J}$ is $3$-generated, i.e., $d\leq 3$.
\qed

\begin{Remark}\rm One suspects that $I$ as above is normal (i.e., $\mathcal{R}(I)$ is a normal domain)
exactly when $d\leq 3\,$; (in other words, exactly when $\mathcal{R}(I)$  is Cohen--Macaulay.)
Of course, in general, a Cohen--Macaulay ideal of codimension $2$ in $k[x,y,z]$ fails very often to even being
integrally closed. The traditional examples are the defining ideals of some of the so-called affine  monomial curves.
For a homogeneous ideal in $k[x,y,z]$ generated in fixed standard degree, one can take the example of \cite{Jar},
namely, $(x^dy^d,x^dz^d,y^dz^d)$ with $d\geq 2$. Although any such ideal is of linear type for  $d\geq 1$,
its integral closure  is not unmixed for $d\geq 2$ -- in any case, by the criterion of \cite{AHA},  this ideal is
not the base ideal of a Cremona map.
This is an entirely different behavior as from
the base ideal of a de Jonqui\`eres map.
Thus, looking from the angle of the defining matrix, the de Jonqui\`eres base ideal is just on the border line, with one column in
standard degree $1$.
\end{Remark}

\subsection{Degree $\leq 4$}

Plane Cremona defined by forms of degree $d\leq 3$ are easily disposed of:

\begin{Proposition}\label{degree3}
Let $I\subset R=k[x,y,z]$ be the base ideal of a plane Cremona map $F\colon
\pp^2\dasharrow \pp^2$ defined by forms of degree $\leq 3$.
Then $F$ is a de Jonqui\`eres map$\,${\rm ;}
in particular, $I$ is saturated.
\end{Proposition}
\demo The case where the degree is $2$ is well-known and completely obvious by the equations of condition
(\ref{eqs_condition}) and B\'ezout theorem. Namely, its weights are $(1,1,1)$.
Therefore, $F$ is a de Jonqui\`eres map, hence the base ideal $I$ is saturated by Corollary~\ref{jonquieres}.

Next let the degree of the three forms be $3$.
This is not even a bit more difficult than the previous case.
We use again the equations of condition and the fact that the number of base points is at least $3$ (\cite[2.6.1]{alberich}).
Since there exists at least  one base point if $d\geq 2$, then for $d=3$ one base point $p$
is such that $\mu_p\geq 2$ (\cite[2.6.8]{alberich}).
Then an immediate calculation yields that the only possible sequence of
multiplicities is $2,1,1,1,1$, in particular there are exactly $5$ proper or infinitely near base points.
Thus, any plane Cremona map of degree $3$ is a de Jonqui\`eres map.
\qed

\medskip

In degree $4$ the situation gets more involved.
Using the equations of conditions and Noether's inequality (\cite[2.6.10]{alberich}),
one has
$$\sum_{p}\mu_p=9,\;\; \sum_{p}\mu_p^2=15,\; \; \mu_1+\mu_2+\mu_3\geq 5,$$
where $\mu_1\geq\cdots \geq\mu_r$.
An elementary calculation yields only two possibilities, namely, either the map is
a de Jonqui\`eres map (of homaloidal type $(4\,;3,1^6)$) or else it is of homaloidal type
$(4\,;2^3,1^3)$.

We know by  Corollary~\ref{degree_at_most_4} that the base ideal of  any Cremona map
in degree $4$ is saturated.
What is left is to describe the invariants of the two alternatives.

The first is dealt with using Corollary~\ref{jonquieres}(a).
For the second possibility, we have the following result:

\begin{Proposition}\label{degree12}
Let $I\subset R=k[x,y,z]$ stand for the base ideal of a Cremona map of
homaloidal type $(4; 2^3,1^3)$, such that the base points are proper and in general position.
Then $I=\mathfrak{F}_{\mathfrak{K}}$ and, in particular, the degree $e(R/I)$ of  $R/I$ is $12$.
\end{Proposition}
\demo
The $3$ quartics defining the map belong to the linear system of all quartics
going through the six points with the prescribed multiplicities $2,2,2,1,1,1$.
The degree of the corresponding $\boldsymbol\mu$-fat ideal $\mathfrak{F}_{\mathfrak{K}}$ is $3+3+3+1+1+1=12$.
 We will use the associativity formula:
$$e(R/I)=\sum_j \ell(R_{P_j}/I_{P_j})\, e(R/P_j),$$
where $P_j$ runs through the minimal primes of $R/I$ and $\ell$ denotes length.
To see how, note that $e(R/P_j)=1$ as $P_j$ is generated by linear forms (since
$k$ is algebraically closed).
Therefore, one has to show $\ell(R_{P_j}/I_{P_j})=3$ if the quartics pass through the
corresponding point with multiplicity $2$, and $\ell(R_{P_j}/I_{P_j})=1$ if the corresponding
point is a simple point on the quartics.

We can assume that $I$ is generated by $3$ quartics each having a double
point at each of the coordinate
points $(0:0:1), (0:1:0), (1:0:0)$, since these points cannot be on a line
(see \cite[Corollary 2.6.9]{alberich}).
Following \cite[Lemma 11.3]{Gibson}, any such quartic can be written in the normal form
$$\lambda x^2y^2+\mu x^2z^2+\nu y^2z^2+2xyz(u_1 x+u_2 y+u_3 z),$$
for suitable values of $\lambda\mu\nu$ not vanishing simultaneously.

Up to another change of coordinates,  we may assume that three quartics generating $I$
are of the form
\begin{eqnarray}\label{quartics}\nonumber
&&x^2y^2+ xyz(t_1 x+t_2 y+t_3 z)=xy(xy+t_1xz+t_2yz+t_3z^2)\\
&& x^2z^2 + xyz(u_1 x+u_2 y+u_3 z)=xz(xz+u_1xy+u_2y^2+u_3yz)\\\nonumber
&&y^2z^2+ xyz(v_1 x+v_2 y+v_3 z)=yz(yz +v_1x^2+v_2 xy +v_3xz),
\end{eqnarray}
for suitable $t_i,u_i,v_i\in k\setminus \{0\}$.
One can further pick special values of the above parameters such that the
parenthetical $2$-forms are reduced, so as to guarantee that the three quartics
admit the remaining $3$ points as simple points and these points be not aligned
(otherwise, the line would be a fixed component thereof).
For example, one may choose
the parameters so as to factor these three $2$-forms into products of linear factors
taken two at the time for each point
(see \cite[2.1.8]{alberich} for a concrete example).

To check $\ell(R_P/I_P)=3$ for a corresponding coordinate point, because of the obvious symmetry, we
can assume that $P=(x,y)$.
Since $z$ is invertible in the localization at $P$, an easy manipulation of the equations in
the format of the right-hand side of (\ref{quartics}) yields $I_P=(x^2,xy,y^2)$,
hence $\ell(R_P/I_P)=3$.
\qed

\subsection{Non-saturated behavior in degree $\leq 7$}

In higher degrees the discussion becomes a lot more involved.
For one thing, there are homaloidal types for which no Cremona map exists
having these types as characteristic (see \cite[5.3.7]{alberich}).
On the bright side, such examples exist only if the number of base points
 exceeds $7$ (\cite[5.3.10]{alberich}).
 Moreover,  for any given proper homaloidal type there exist simple Cremona maps
 whose base points can be generically chosen (\cite[5.3.5]{alberich}).
For this reason, we will focus only on simple Cremona maps.

As explained in Subsection~\ref{analogues}, the multiplicity (geometric degree) of the base locus of a simple Cremona map
with weighted cluster $\mathcal{K}=(K,\boldsymbol\mu)$ is
bounded below by the degree of the associated $\boldsymbol\mu$-fat ideal $\mathfrak{F}_{\mathcal{K}}$.
We will say that a simple Cremona map has {\em minimum multiplicity}
if the multiplicity of its base locus attains its lower bound.
Clearly, a similar notion could be introduced for non-simple maps by resorting to the
multiplicity of the associated full ideal of curves $I_{\mathcal{K}}$, but, unfortunately, we know no closed formula
for the latter.

\begin{Proposition}\label{pminimaldegree}
The homaloidal type of a  simple Cremona map of minimum multiplicity is one of the following:
\begin{equation}\label{minimal_multiplicity_types}
(5\,;2,2,2,2,2,2), \;(4\,;2,2,2,1,1,1),\; (3\,;2,1,1,1,1), \;(2\,;1,1,1)
\end{equation}
\end{Proposition}
\demo Write $\boldsymbol\mu=\{\mu_1,\cdots,\mu_n\}$.
By the associativity formula, the degree of the base locus is
$$e(R/I)=\sum_{j=1}^n \ell(R_{P_j}/I_{P_j})\, e(R/P_j)=\sum_{j=1}^n \ell(R_{P_j}/I_{P_j}),$$
since the base field is assumed to be algebraically closed.
The minimal multiplicity assumption means that $\sum_j \ell(R_{P_j}/I_{P_j})= \sum_j \ell((R/P_{j}^{\mu _j})_{P_j}) $,
where $1\leq j\leq n$. But the minimal number of generators of $I_{P_j}$ is $3$ whereas $(P_{j}^{\mu _j})_{P_j}=({P_j}_{P_j})^{\mu_j}$
 requires at least $4$ minimal  generators if $\mu_j\geq 3$.
  Therefore, one must have $\mu_j\leq 2$ for $1\leq j\leq n$.

Thus, up to ordering, we may assume that $\{\mu_1,\cdots,\mu_n\}=
\{\underbrace{1,\ldots,1}_j; \underbrace{2,\ldots,2}_k\}$. Then the equations of condition read
\begin{equation}\label{eqsystem}
\left\lbrace
           \begin{array}{c l}
               j+2k=3d-3\\
              j+4k=d^2-1.
           \end{array}
         \right.
\end{equation}
This system gives $j=(5-d)(d-1)$ which yields $d\leq 5$ as $j\geq 0$.
The asserted alternatives now follow by an easy case by case calculation.
\qed

\begin{Theorem}\label{tcharactrization} Let $F$ be a simple plane Cremona map of degree $d\leq 7$ with base ideal $I$.
Assume that $I$ is not saturated {\rm (}hence, $d\geq 5${\rm )}.
\begin{itemize}
\item[\rm (i)] If $d=5$, $F$ has homaloidal type $(5\,;2^6):=(5\,;2,2,2,2,2,2)$, $e(R/I)=18$, and the minimal resolution of $R/I$
has the form
$$0\rightarrow R(-9)\rightarrow R^3(-8)\rightarrow R^3(-5)\rightarrow R.$$
\item[\rm (ii)] If $d=6$, $F$ has homaloidal type $(6\,;4,2^4,1^3)$, $e(R/I)=26$, and the minimal resolution of $R/I$
has the form
$$0\rightarrow R(-11)\rightarrow R^2(-10)\oplus R(-9) \rightarrow R^3(-6)\rightarrow R.$$
\item[\rm (iii)] If $d=7$, $F$ has homaloidal type $(7\,;5,2^5,1^3)$ or else $(7\,;4,3^2,2^3,1^2)$. In the second case
 $e(R/I)=36$ and the minimal resolution of $R/I$ has one of the following forms:
    $$0\rightarrow R(-12)\rightarrow R^3(-11) \rightarrow R^3(-7)\rightarrow R$$
 or
 $$0\rightarrow R(-13)\rightarrow R^2(-12)\oplus R(-10) \rightarrow R^3(-7)\rightarrow R.$$
\end{itemize}
\end{Theorem}
\demo By the previous sections, one has only to consider the range $5\leq d\leq 7$.

(i) According to Proposition~\ref{P567}(i) the minimal free resolution of the base ideal of a plane Cremona map (not necessarily simple)
 of degree $5$ has the form mentioned there. Thus, for all such maps $e(R/I)=18$. This is the minimum value for simple maps
 (Proposition~\ref{cupperbound}). The result, now, follows from Proposition~\ref{pminimaldegree}.

(ii) By Proposition~\ref{P567}(ii) there are two virtual resolutions of $R/I$. By Proposition~\ref{pminimaldegree},
the first of these alternatives conflicts with the content of Proposition~\ref{P567}(ii) as applied to simple plane maps.
The second alternative says that $e(R/I)=26$, which is exactly one more than the minimum multiplicity.
Thus among the inequalities $\ell((R/I)_{P_i})\geq \ell((R/P_{i}^{\mu _i})_{P_i})$ for $i=1,\cdots,n$, only one of them is strict.
Hence, by the same argument as in the proof of Proposition~\ref{pminimaldegree}, we may assume that $\mu_1\geq 3$
and $\mu_i\leq 2$ for all $i\geq 2$.
According to Corollary~\ref{jonquieres}(a),  $\mu_1\leq 4$. Thus, $\mu_1\in\{3,4\}$. If $\mu_1=3$, the equation of conditions lead to
the following system (with the similar notations as in (\ref{eqsystem}))
$$\left\lbrace
           \begin{array}{c l}
               j+2k=3d-3-3=12\\
              j+4k=d^2-1-9=26
           \end{array}
         \right.
$$
This system has no solution in nonnegative integers. For $\mu_1=4$, the solution yields the homaloidal type $(6\,;4,2^4,1^3)$.

(iii) The case where $d=7$ is more involved. Again Proposition~\ref{P567}(iii) implies that $33\leq e(R/I)\leq 36$, where $33$ is
the minimum multiplicity. Then among the inequalities $\ell((R/I)_{P_i})\geq \ell((R/P_{i}^{\mu _i})_{P_i})$ for $i=1,\cdots,n$
at most three ones are strict. Therefore, we may assume that $3\leq \mu_1\leq 5$ and $\mu_i\leq 2$ for $i\geq 4$ as in the proof
of Proposition~\ref{pminimaldegree} and by Corollary~\ref{jonquieres}(a).
We now deal separately with the three alternatives for $\mu_1$.

If $\mu_1=3$ we have $(d-\mu_1)/2=(7-3)/2=2$, hence $\mu_2=\mu_3=3$ (\cite[Definition 8.2.1 and Lemma 8.2.6]{alberich}).
Plugging these data into the equations of condition yields no solution in positive integers.

Let $\mu_1=5$. In this case the equations of condition yield no solution in positive integers for $4\leq \mu_2\leq 5$ or $\mu_2=\mu_3=3$.
Then the triple $(\mu_1,\mu_2,\mu_3)\in \{(5,3,2),(5,2,2)\}$. Plugging the first of these into the equations
of condition yields the homaloidal type $(7\,;5,3,2^2,1^6)$. In a similar vein, $(5,2,2)$ yields the homaloidal type $(7\,;5,2^5,1^3)$.
Now, an application of the Hudson's test shows that the first of these homaloidal types is improper.
Therefore, we are left with $(7\,;5,2^5,1^3)$ as the only possibility.

Finally, let $\mu_1=4$. The triple $(\mu_1,\mu_2,\mu_3)$ belongs to the following list
$$(4,4,4),(4,4,3),(4,4,2),(4,3,3),(4,3,2),(4,2,2).$$
Now, let the set $\{\mu_i:i\geq4\}$ have $j$ elements equal to $1$
 and $k$ elements equal to $2$. Consider the system given by the equations of condition:
$$\left\lbrace
           \begin{array}{c l}
               j+2k=3d-3-(\mu_1+\mu_2+\mu_3)\\
              j+4k=d^2-1-(\mu_1^2+\mu_2^2+\mu_3^2)
           \end{array}
         \right.
$$
The solution pair $(j,k)$ is, respectively, $(-3,12),(0,7),(3,4),(3,2),(5,-1),(7,-4)$.
Clearly, the first and the last two pairs are absurd.
Next applying Hudson's test  to the first two possible homaloidal types shows that they are improper.
Hence, the remaining possibility is $(7\,;4,3^2,2^3,1^2)$.

\medskip

 As to the form of the resolution of $R/I$ is, notice that, for this homaloidal type, the first three multiplicities exceed $2$.
It follows that $e(R/I)\geq 3 +33=36$, which is the maximal possible value.
Then among the virtual resolutions in Proposition~\ref{P567}(iii) only the last two satisfy $e(R/I)=36$, as claimed.
\qed

\medskip

We note that the converse of item (ii) in the above theorem is not true in general, at least if one does not assume that $F$ is simple.
\begin{Example}\rm The following non-simple Cremona map appears in \cite[Example 2.1.14]{alberich}:
$$\bigl((x^3 - yz (y + x)) (x^2 - yz) (y + x)\,:\,(x^2- yz) x^2 (x + y)^2\,:\,x^3 (x^3 -yz(x + y))\bigr).$$
\end{Example}
The homaloidal type of this map is
$(6\,; 4,2^4,1^3)$, but it has a saturated base ideal.


\subsection{Structure in degree $5$}

Degree $5$ has special arithmetic features, such as being the only degree $d\geq 2$ having the property that
\begin{equation}
\sum_{p}\mu_p=1/2 \sum_{p}\mu_p^2.
\end{equation}
Further, an elementary scrutiny in the equations of condition yields that the only homaloidal types are
\begin{enumerate}
\item[{\rm (a)}] $4,\underbrace{1,\ldots,1}_8$ (de Jonqui\`eres type)
\item[{\rm (b)}] $3,3,\underbrace{1,\ldots,1}_6$
\item[{\rm (c)}] $3,2,2,2,1,1,1$
\item[{\rm (d)}] $\underbrace{2,\ldots,2}_6$ (symmetric type)
\end{enumerate}

Of these, only type (b) is not proper (cf. \cite[Example 5.3.7]{alberich})
and the classical result the fact that for $7$ points or less every homaloidal type is proper).
Some of this can be readily checked by Hudson's test.

Our purpose is to give a precise characterization of the three proper homaloidal types
in terms of the homological features of the corresponding base ideals.

\begin{Theorem}\label{proper_types_degree5}
Let   $F\colon\pp^2\dasharrow \pp^2$ stand for a simple plane Cremona map of degree $5$ whose base points
are in general position, and let $I\subset R=k[x,y,z]$ denote its base ideal.
Then
\begin{enumerate}
\item[{\rm (i)}] $F$ is a de Jonqui\`eres map if and only if  $R/I$ is Cohen--Macaulay of degree  $21$
with resolution of the form
\begin{equation}\label{res_of_jonquieres_deg5}
0\rar R(-9)\oplus R(-6)\lar R(-5)^3\lar R\,{\rm ;}
\end{equation}
\item[{\rm (ii)}] $F$ has type $(5\,;3,2^3,1^3)$ if and only if  $R/I$ is Cohen--Macaulay of degree $19$
with resolution of the form
\begin{equation}\label{res_of_extratype_deg5}
0\rar R(-8)\oplus R(-7)\lar R(-5)^3\lar R\,{\rm ;}
\end{equation}
\item[{\rm (iii)}] $F$ is symmetric of type $(5\,;2^6)$ if and only if $R/I$ is non Cohen--Macaulay of
degree $18$ {\rm (m}inimal possible{\rm )} with resolution of the form
\begin{equation}\label{res_of_symmetric_deg5}
0\rar R(-9)\stackrel{\psi_2}{\lar} R(-8)^3\stackrel{\psi_1}{\lar} R(-5)^3\lar R.
\end{equation}
{\sc Supplement.} The maps $\psi_1,\, \psi_2$ come out of the minimal free  resolution of the corresponding fat ideal
$\mathfrak{F}$
$$0\rar R(-7)^3\stackrel{\phi}{\lar} R(-5)^3\oplus R(-6)\lar R$$
as follows: up to shifts, $\psi_1$ is the the restriction of the dual map $\phi^t$, and $\psi_2$ is the tail of the
Koszul complex on three $k$-linearly independent linear forms generating the coordinates of the syzygies of the
additional minimal generator of $\mathfrak{F}$ of degree $6$.
\end{enumerate}
\end{Theorem}
\demo
We note at the outset that there are three virtual resolutions in this degree.
Indeed, if $I$ is not saturated then according to Theorem~\ref{tcharactrization} the minimal resolution must be of the form
$$0\rar R(-9)\lar R(-8)^3\lar R(-5)^3\lar R.$$
If $I$ is saturated then one easily sees that the only possible minimal resolutions are
\begin{equation}\label{res5sat1}0\rightarrow R(-9)\bigoplus R(-6)\rightarrow R^3(-5)\rightarrow R\end{equation} or
\begin{equation}\label{res5sat2}0\rightarrow R(-8)\bigoplus R(-7)\rightarrow R^3(-5)\rightarrow R\end{equation}

Let us first prove (iii).
Since we are assuming that the Cremona map is simple, the ``if" part  follows from Proposition~\ref{pminimaldegree}
(or also from Theorem~\ref{tcharactrization}).

We now show the ``only if'' part of (iii).
Thus assume that the homaloidal type is  $(5\,;2^6)$.
Let $\mathfrak{F}\subset k[x,y,z]$ stand for the associated $\mathbf{2}$-fat ideal  based on the given points.

The homological nature of $\mathfrak{F}$  has been determined in \cite{Fich}:
\begin{equation}\label{six2_fat_points}
0\rar R(-7)^3\stackrel{\phi}{\lar} R(-5)^3\oplus R(-6)\lar R\rar R/\mathfrak{F}.
\end{equation}
In particular, $\mathfrak{F}=(I,f)$, where $f$ is a $6$-form and there are three linear forms $\ell_1,\ell_2,\ell_3$
such that $\ell_i\,f\in I$ for every $i$.

Suppose that $\ell_1,\ell_2,\ell_3$ are $k$-linearly dependent.
Then a suitable $k$-linear combination of syzygies gives a syzygy of the
generators of $I$ in standard degree $2$.
Such a syzygy would necessarily be of minimal degree among the syzygies of $I$ since the latter
cannot have linear syzygies because the resolution (\ref{six2_fat_points}) is minimal.
But, among the three virtual resolutions above, the only one having such a minimal syzygy is (\ref{res5sat2}).
We now argue that this is impossible.

Indeed, let $\{p_1,\ldots,p_n\}$ denote the set of proper base points and  $\{P_1,\ldots,P_n\}$ the
corresponding defining prime ideals. Now, for this free resolution, $e(R/I)=19$ which is exactly one more than the minimum degree.
Therefore, there is exactly one prime ideal, say $P$, among $P_1,\ldots,P_n$ such that $\ell((R/I)_P)=\ell((R/P^2)_P)+1$.
Without loss of generality we may assume that $P=(x,y)$. It then follows that $\ell((P^2/I)_P)=1$,
that is, $(P^2/I)_P\simeq (R/P)_P$; in particular,  $P^3_P\subset I_P$.
Let us consider the three cases $(P^2/I)_P=(\overline{xy}/1)$, $(P^2/I)_P=(\bar{x}^2/1)$ or
$(P^2/I)_P=(\bar{y}^2/1)$. Say, $(P^2/I)_P=(\overline{xy}/1)$.
Then, $(\bar{x}^2/1)=\frac{a_1(z)}{a_2(z)}(\overline{xy}/1)$ and
$(\bar{y}^2/1)=\frac{b_1(z)}{b_2(z)}(\overline{xy}/1)$ where $a_i(z)$ and $b_i(z)$
are polynomials in $k[z]$. Solving these equations, one finds polynomials $c_i(z)$ and $d_i(z)$ for $i=1,2$
such that $c_1(z)x^2-c_2(z)xy\in I$ and $d_1(z)y^2-d_2(z)xy\in I$. The fact that $I$ is a homogeneous ideal then
implies that there are integers $n,m$ and elements $b,b' \in k$ such that $g_1=z^n(x^2-bxy)$ and $g_2=z^m(y^2-b'xy)$
belong to $I$. Thence, all of the base points of $F$ are contained in the variety defined by $g_1$ and $g_2$
which consists of two lines, (provided $bb'=1$). This, however, contradicts the assumption that the base points
are in general position.
The remaining two other cases are dealt with in an entirely  similar fashion.

\smallskip

Therefore, it must be the case that $\ell_1,\ell_2,\ell_3$ are $k$-linearly independent, hence
generate the maximal ideal $(x,y,z)$.
It follows that,  $f$ drives $(x,y,z)$ inside $I$.
This means that depth$(R/I)=0$ as $f\not\in I$.
In other words, $R/I$ is not Cohen--Macaulay and $I^{\rm sat}= \mathfrak{F}$; in particular, $e(R/I)=18$, as claimed in
the statement.
This concludes the proof of (iii).

\smallskip

The ``only if'' implication of (i) was already shown in Corollary~\ref{jonquieres}.

We now prove the ``only if'' implication of (ii).
The free resolution of the associated fat ideal $\mathfrak{F}$ is known (\cite{Harb}):
\begin{equation}\label{fat_extratype_deg5}
0\rar R(-7)^3\stackrel{\phi}\lar R(-5)^3\oplus R(-6)\lar R.
\end{equation}

Therefore, up to a projective coordinate change, we may assume that $\mathfrak{F}$ is
is defined by the three minimal generators of degree $5$ of $I$ and one extra generator of degree $6$, say $f$.
Focusing on the entries of $\phi$, we see that there are $3$ linear forms $\ell_1,\ell_2,\ell_3$
such that $\ell_i\,f\in I$ for every $i$.
By Theorem~\ref{tcharactrization} in this type, $(5\,;3,2^3,1^3)$,  $I$ is saturated. This implies that
$\ell_1,\ell_2,\ell_3$ are $k$-linearly dependent. Thus, by a  projective coordinate change,
one may assume that $\ell_3=0$. Therefore the third column in $\phi$ is a syzygy of (standard) degree $2$ of
the generators of $I$.
Since $R/I$ is Cohen--Macaulay and $I$ is generated in degree $5$ the resolution has got to be of the form
indicated in (\ref{res_of_extratype_deg5}).
From this we easily compute the degree of $R/I$, as stated.

\smallskip

In order to complete the proof it now suffices to prove the ``if'' parts of  (i) and (ii).
But since we have completed the proofs of the ``only if'' parts of all three items,
this follows automatically.

\smallskip

It remains to deal with the supplement of item (iii), to explain the nature of the maps in the resolution of the non-saturated case.
First dualize (\ref{six2_fat_points}) into $R$:
$$
0\rar J^*\simeq R\lar R(5)^3\oplus R(6)\stackrel{\phi^t}{\lar} R(7)^3.
$$
Consider the restriction $\psi:=\phi^t\restr _{R(5)^3}: R(5)^3\lar R(7)^3$.
Then consider the left tail of the Koszul  complex on
the regular sequence $\ell_1,\ell_2,\ell_3$ shifted by $5$:
$$0\rar R(3)\lar R(4)^3\lar R(5)^3$$
Call $\psi$ the composite of the right most map with the restriction $\phi^t\restr _{R(5)^3}: R(5)^3\lar R(7)^3$
of $\phi^t$ to $R(5)^3$.
The following sequence of maps obtains:
$$0\rar R(3)\lar R(4)^3\stackrel{\psi}{\lar} R(7)^3.$$
Now shift by $-12$ to get
$$0\rar R(-9)\lar R(-8)^3\stackrel{\Psi}{\lar} R(-5)^3,$$
where $\Psi=\psi(-12)$ for lighter reading.
We claim this resolves $I$.
Proving this involves  first showing that it is a complex.
For it, the only missing piece is that $\Psi(R(-8)^3)$ are syzygies of the three quintics generating $I$.
But this is clear since the $i$th row of $\psi$ is
$$
(h_1^{(i)}, \, h_2^{(i)},\, h_3^{(i)})
\left(
\begin{array}{ccc}
0&\ell_3 & -\ell_2\\
-\ell_3& 0 & \ell_1\\
\ell_2 & -\ell_1 & 0
\end{array}
\right)=
(-h_2^{(i)}\ell_3 + h_3^{(i)}\ell_2,\, h_1^{(i)}\ell_3-h_3^{(i)}\ell_1,\, -h_1^{(i)}\ell_2+h_2^{(i)}\ell_1),
$$
where $(h_1^{(i)}, \, h_2^{(i)},\, h_3^{(i)})$ denotes the $i$th row of $\phi^t$ chopping off the fourth coordinate,
But the right most row above is obtained from the three syzygies of $J$ by multiplying by $\ell_j$ and by $\ell_k$
the syzygy involving $\ell_m$ as coordinate, where $\{j,k.m\}=\{1,2,3\}$, then subtracting in the obvious way.

To use the acyclicity criterion as formulated by Buchsbaum--Eisenbud (\cite[Theorem 20.9]{E})
it remains to argue that $I_2(\Psi)$ has height $\geq 2$.
But $\Psi$ is the transpose of the composite
$$R(-7)^3\stackrel{\phi}{\lar} R(-5)^3\oplus R(-6)\lar R(-5)^3,$$
where the rightmost map is projection onto the first factor.
Therefore, one has to show that the $2$-minors of the first $3$ rows of the latter matrix has codimension $\geq 2$.
But this is clear as this ideal contains the ideal of $3$-minors of $\phi$, which has codimension $2$.
This shows that the complex is exact.

\qed

\begin{Corollary}\label{bi_homaloidness}\rm
A Cremona map of degree $\leq 5$ and its inverse have the same homaloidal type.
\end{Corollary}
\demo
Since there are only three homaloidal types and the inverse of a de Jonqui\`eres map is
a de Jonqui\`eres map, it suffices to argue that a Cremona map  and its inverse
cannot be of type $(5\,;3^2,1^6)$ and $(5\,;2^6)$,
respectively.
But this follows from Clebsch theorem \cite[Theorem 3.3.2]{alberich}.
\qed

\medskip

It could be of some interest to give an independent proof of the preceding corollary,
staged in the spirit of the rest of the paper.
By the results obtained so far, it would suffice to know that a Cremona map of degree $5$ and its inverse
have base ideals of the same (scheme) degree.
A parallel challenge is to decide when the saturation $I^{\rm sat}$ of the base ideal $I$ of a simple Cremona map $F$
coincides with the associated  $\boldsymbol\mu$-fat ideal (recall that, from Corollary~\ref{sat_is_closure},
this implies that $I^{\rm sat}$ is the integral closure of $I$).
Unfortunately,  many Cremona maps do not share any of the two  properties,
which is already the case in degree $6$ \cite[Example 2.1.14]{alberich}.

\begin{Remark}\rm In principle it ought to be possible to explicitly enumerate all
homaloidal nets stemming from ideals of fat points up to any given upper bound
for the multiplicities of the given set of proper points.
For any $3$ linearly independent forms of a fixed degree thus found, one could
apply Hudson's test.
It would be interesting to understand the complexity of one such algorithm.
\end{Remark}

\subsection{An algebraic test}

One often needs a more algebraic criterion to decide when a good candidate
is indeed a Cremona map.
The following statement gives such an alternative test for a special class of plane rational maps.
It could be used to prove that the linear system of quartics through six  points
with multiplicities $2,2,2,1,1,1$ defines a Cremona map, regardless of the point configuration.

\begin{Proposition}\label{ideal_squared_replaces_telescopic}
Let $J\subset R$ be a homogeneous ideal having a minimal free resolution of the form
$$0\rar R(-6)^2\stackrel{\phi}{\lar} R(-4)^3\lar R$$
and analytic spread $3$.
If $R/J^2$ is also Cohen--Macaulay then the generators of $J$ define a Cremona map of
$\pp^2$ of degree $4$.
\end{Proposition}
\demo
Let $\phi$ also denote a matrix
of the map $\phi$. Say,
$$\varphi=\left(\begin{array}{rr}
\alpha_1 & \alpha_2 \\
\beta_1 & \beta_2 \\
\gamma_1 & \gamma_2\\
\end{array} \right)
$$
It produces trivially the following matrix of syzygies of
$J^2$:
$$\phi^{[2]}=\left(\begin{array}{rrrrrr}
\alpha_1 & \alpha_2 & 0 & 0 & 0 & 0\\
\beta_1 & \beta_2 & \alpha_1 & \alpha_2 & 0 & 0\\
\gamma_1 & \gamma_2 & 0 & 0 & \alpha_1 & \alpha_2\\
0 & 0 & \beta_1 & \beta_2 & 0 & 0\\
0 & 0 & \gamma_1 & \gamma_2 & \beta_1 & \beta_2\\
0 & 0 & 0 & 0 & \gamma_1 & \gamma_2
\end{array} \right)
$$
We claim that $\phi^{[2]}$ has maximal possible rank, namely,
$5$. The determinant of the first $5$ columns and last $5$ rows is
$\pm \gamma_1(\beta_1\gamma_2-\beta_2\gamma_1)^2$. By symmetry, we
find a similar determinant whose value is $\pm
\gamma_2(\beta_1\gamma_2-\beta_2\gamma_1)^2$. Since $\gamma_1,
\gamma_2$ cannot simultaneously vanish, we are through.

Now, since $\rk \phi^{[2]}=5$ and the total syzygy matrix $\Phi$ of
$J^2$ is $5\times 6$ of rank $5$ - because by assumption the analytic spread of $J$ is $3$
and $R/J^2$ is Cohen--Macaulay --
then the columns of $\Phi$ have
standard degree $\leq 2$. However, its $5\times 5$ minors are the
generators of $J^2$  which are of degree
$8$, the only possibility that adds up correctly is that the columns
of $\Phi$ be of degrees $1,1,2,2,2$.

This argument shows that the ideal $J^2$ has
two independent syzygies with linear coordinates.
Now consider a presentation of the Rees algebra $\mathcal{R}_R(J)\simeq R[t,u,v]/\mathcal{J}$,
with $\mathcal{J}$ a bihomogeneous ideal in the standard bigrading of $R[t,u,v]=k[x,y,z,t,u,v]$.
Then the two linear syzygies of $J^2$ induce generators of $\mathcal{J}$ of bidegree
$(1,2)$ generating a subideal of codimension $2$.
We can now apply \cite[Proposition 3.9]{AHA}.
\qed

\medskip

It seems natural, in the present stringent setup, to pose:

\begin{Question}\rm
Under the hypotheses of Proposition~\ref{ideal_squared_replaces_telescopic} is the converse true,
i.e., if $J$ defines a Cremona map must $R/J^2$ be Cohen--Macaulay?
\end{Question}
Observe that the assumption on $J^2$ is delicate. Thus, the general $3\times 2$ matrix $\phi$  with $2$-forms
as entries defines a Cohen--Macaulay ideal $J$ whose square is not saturated and has minimal resolution
of the form $0\rar R(-12)\rar R(-10)^6\rar R(-8)^6\rar R$.
As a matter of fact, the syzygies of any power $J^m$, with $m\geq 2$, have degree at least $2$,
hence the rational map defined by a set of minimal generators of $J$ is not Cremona by the criterion of
\cite{AHA}.
There is an easy background explanation for the map failing to be Cremona and that is the fact
that $J$ is a radical ideal, hence it is the fat ideal of points with multiplicities unit -- this
defines a Cremona map only in degrees $\leq 2$.

\begin{Question}\rm
If $\phi$ is not general, an interesting side question is what are the homological properties of $J$
that trigger the existence of a unique linear syzygy among all powers of $J$? By \cite{AHA} such an ideal does not define
a Cremona map, however the nature of the base points is subtler.
\end{Question}


\noindent {\bf Authors' addresses:}

\medskip

{\em S.H.Hassanzadeh},  Departamento de Matem\'atica, CCEN, Universidade Federal
de Pernambuco,
Cidade Universit\'aria, 50740-540 Recife, PE, Brazil {\em and}
  Faculty of Mathematical Sciences and Computer, Tarbiat
Moallem University, Tehran, Iran .\\ email: hamid@dmat.ufpe.br\\

{\em A.Simis},  Departamento de Matem\'atica, CCEN, Universidade Federal
de Pernambuco,
Cidade Universit\'aria, 50740-540 Recife, PE, Brazil.\\ email: aron@dmat.ufpe.br

\end{document}